\documentclass[11pt]{amsart}

\usepackage{amsfonts}
\usepackage{amssymb}
\usepackage{amscd}

\theoremstyle{plain}

\newtheorem{theorem}{Theorem}[section]
\newtheorem{proposition}[theorem]{Proposition}
\newtheorem{corollary}[theorem]{Corollary}
\newtheorem{lemma}[theorem]{Lemma}
\theoremstyle{definition}
\newtheorem{definition}[theorem]{Definition}
\newtheorem{remark}[theorem]{Remark}

\newtheorem{conjecture}[theorem]{Conjecture}
\newtheorem{conjecture/question}[theorem]{Conjecture/Question}

\newtheorem{remark/definition}[theorem]{Remark/Definition}
\newtheorem{terminology/notation}[theorem]{Terminology/Notation}

\def\rmapdown#1{\Big\downarrow
   \rlap{$\vcenter{\hbox{$\scriptstyle#1$}}$ }}
 \theoremstyle{remark}

\begin{document}

\title{\bf Gaussian maps, Gieseker-Petri loci and large theta-characteristics}

\author[G. Farkas]{Gavril Farkas}
\address{Department of Mathematics, Princeton University,
Fine Hall, Princeton, NJ 08544}
\email{{\tt gfarkas@math.princeton.edu}}
\thanks{Research partially supported by the NSF Grant DMS-0140520}

%\date{\today}
\maketitle

\setlength{\oddsidemargin}{0 in}
\setlength{\footskip}{.3 in}
\setlength{\headheight}{.3 in}
\setlength{\textheight}{8.5 in}
\setlength{\parskip}{.1 in}

\newcommand{\marginlabel}[1]%
  {\mbox{}\marginpar{\raggedleft\hspace{0pt}\bfseries\sf#1}}

\def\ZZ{{\mathbb Z}}
\def\GG{{\mathbb G}}
\def\QQ{{\mathbb Q}}
\def\PP{{\textbf P}}
\def\OO{{\mathcal O}}
\def\cX{\mathcal{X}}
\def\cC{\mathcal{C}}
\def\cA{\mathcal{A}}
\def\F{\mathcal{F}}
\def\cE{\mathcal{E}}
\def\G{\mathcal{G}}
\def\K{\mathcal{K}}
\def\L{\mathcal{L}}
\def\S{\mathcal{S}}
\def\J{\mathcal{J}}
\def\cM{\mathcal{M}}
\def\cZ{\mathcal{Z}}
\def\cU{\mathcal{U}}
\def\cQ{\mathcal{Q}}
\def\cI{\mathcal{I}}
\def\Pic0{{\rm Pic}^0(X)}
\def\pxi{P_{\xi}}
\def\ff{\overline{\mathcal{F}}}
\def\mm{\overline{\mathcal{M}}}
\def\ss{\overline{\mathcal{S}}}

\markboth{G. FARKAS}
{Gaussian maps and theta-characteristics}

\section{Introduction}
For an integer $g\geq 1$ we consider the moduli space $\S_g$ of
smooth spin curves parametrizing pairs $(C,L)$, where $C$ is a
smooth curve of genus $g$ and $L$ is a theta-characteristic, that
is, a line bundle on $C$ such that $L^2\cong K_C$. It has been
known classically that the natural map $\pi: \S_g \rightarrow
\cM_g$ is finite of degree $2^{2g}$ and that $\S_g$ is a disjoint
union of two components $\S_g^{even}$ and $\S_g^{odd}$
corresponding to even and odd theta-characteristics. A
geometrically meaningful compactification $\ss_g$ of $\S_g$ has
been constructed by Cornalba by means of stable spin curves of
genus $g$ (cf. \cite {C}).  The space $\ss_g$ and more generally
the moduli spaces $\ss_{g,n}^{1/r}$ of stable $n$-pointed $r$-spin
curves of genus $g$, parametrizing pointed curves with $r$-roots
of the canonical bundle, have attracted a lot of attention in
recent years, partly due to a conjecture of Witten relating
intersection theory on $\ss_{g,n}^{1/r}$ to generalized KdV
hierarchies (see e.g. \cite{JKV}).
\newline
\indent For each $g,r\geq 0$ one can define the locus
$$\S_g^r:=\{(C,L)\in \S_g: h^0(L)\geq r+1 \mbox{ and } h^0(L)\equiv r+1 \mbox{ mod 2}\}.$$
We also set $\cM_g^r:=\pi(\S_g^r)$. It has been proved by  Harris
that each component of $\S_g^r$ has dimension $\geq 3g-3-{r+1\choose 2}$
(cf. \cite{H}). This bound is known to be sharp when $r$ is very
small: it is a classical result that $\S_g^1$ is a divisor in
$\S_g$, while for $r=2,3$ we have that $\S_g^r$ has pure
codimension $r(r+1)/2$ in $\S_g$ for all $g\geq 8$ (cf.
\cite{T1}). On the other hand clearly the bound is far from
optimal when $r$ is relatively large with respect to $g$ in the
sense that there are examples when $\S_g^r\neq \emptyset$ although
$3g-3- {r+1 \choose 2}$ is very negative: the hyperelliptic locus
$\mathcal{H}_g\subset \cM_g$ is contained in $\cM_{g}^{[(g-1)/2]}$
and there are Castelnuovo extremal curves $C\subset \PP^r$ of
genus $3r$ such that $K_C=\OO_C(2)$, which gives that
$\S_{3r}^r\neq \emptyset $ for all $r\geq 3$ (see e.g.
\cite{CdC}). It is thus natural to ask to what extent Harris'
bound is sharp. We give a partial answer to this question by
proving the following:

\begin{theorem}\label{th1}
For $1\leq r\leq 11, r\neq 10$, there exists an explicit integer
$g(r)$ such that for all $g\geq g(r)$ the moduli space $\S_g^r$
has at least one component of codimension ${r+1 \choose 2}$ in $\S_g$.
The general point $[C,L]$ of such a component corresponds to a
smooth curve  $C\subset \PP^r$, with $L=\OO_C(1)$ and $K_C=\OO_C(2)$.
\end{theorem}

For a precise formula for $g(r)$ we refer to Section 3. We
conjecture the existence of a component of $\S_g^r$ of codimension
${r+1 \choose 2}$ for any $r\geq 1$ and $g\geq {r+2\choose 2}$ and we
indicate a way to construct such a component (see Conjecture \ref{white}). Theorem
\ref{th1} is proved inductively using the following result:

\begin{theorem}\label{th2}
We fix integers $r,g_0\geq 1$. If $\S_{g_0}^r$ has a component of
codimension ${r+1 \choose 2}$ in $\S_{g_0}$, then for every $g\geq
g_0$, the space
 $\S_g^r$ has a component of codimension ${r+1 \choose 2}$ in $\S_g$.
\end{theorem}

To apply Theorem \ref{th2} however, one must have a starting case
for the inductive argument. This is achieved  by carrying out an
infinitesimal study of the loci $\S_g^r$ which will relate
theta-characteristics to Gaussian maps on curves.  Recall that for
a smooth curve $C$ and a line bundle $L$ on $C$, the
\emph{Gaussian} or \emph{Wahl} map $\psi_L:\wedge ^2 H^0(L)\rightarrow
H^0(K_C\otimes L^2)$ is defined essentially by
$$\psi_L(s\wedge t):=s\ dt-t\ ds.$$
The map $\psi_L$ has attracted considerable interest being
studied especially in the context of deformation theory (see
\cite{W1}  and the references therein).  Wahl proved
the remarkable fact that if $C$ sits on a $K3$ surface then
$\psi_{K_C}$ cannot be surjective, which should be contrasted with
the result of Ciliberto, Harris and Miranda saying that
$\psi_{K_C}$ is surjective for the general curve $C$ of genus
$g=10$ or $g\geq 12$ (cf. \cite{CHM}). In a completely different direction,
in a previous work we made essential use of the Gaussian map
$\psi_{K_C}$ for $g=10$ to construct a counterexample to the
Harris-Morrison Slope Conjecture on effective divisors on $\mm_g$
(cf. \cite{FP}).

There are several powerful criteria in the literature ensuring the
surjectivity of $\psi_L$ when $L$ has large degree (see e.g.
\cite{Pa}, Theorem G), but very little seems to be known about
when is the map $\psi_L$ injective, or more generally, what is the
behaviour of  $\psi_L$ when the line bundle $L$ is special (cf.
Question 5.8.1 in \cite{W1}). In Section 5 we go some way towards
answering this question by showing the following:
\begin{theorem}\label{th3}
For the general curve $C$ of genus $g$ and for any line bundle $L$
on $C$ of degree $d\leq g+2$, the Gaussian map $\psi_L$ is
injective.
\end{theorem}

We refer to Theorem \ref{th4}  for a more general statement that
bounds the dimension of $\mbox{Ker}(\psi_L)$ even when $d>g+2$. In
the case when $L$ is a very ample line bundle giving an embedding
$C\subset \PP^n$, Theorem \ref{th3} can be interpreted as saying
that the associated curve $C\rightarrow \PP^N$ obtained by
composing the Gauss map $C\rightarrow G(2,n+1)$, $C\ni p\mapsto
\mathbb T_p(C)$, with the Pl\"{u}cker embedding of the
Grassmannian of lines, is nondegenerate. Alternatively one can
read this result in terms of (absence of) certain
self-correspondences on the general curve $C$ (see Proposition
\ref{corr}).

In Section 4 we relate the Gieseker-Petri loci on $\mm_g$ to the
moduli spaces $\S_{g,n}^r$ of $n$-pointed spin curves consisting
of collections $(C,p_1,\ldots, p_n, L)$, where $(C, p_1,\ldots
p_n)\in \cM_{g,n}$ and $L$ is a degree $k$ line bundle on $C$ such
that $L^2\otimes \OO_C(p_1+\cdots+p_n)=K_C$ and $h^0(L)\geq r+1$.
Here of course we assume that $2k+n=2g-2$.

We recall that the Giseker-Petri Theorem asserts that for a
general curve $C$ of genus $g$ and for any line bundle $L$ on $C$,
the  map $\mu_0(L):H^0(L)\otimes H^0(K_C\otimes L^{-1})\rightarrow
H^0(K_C)$ is injective (see e.g. \cite{EH2}). It is
straightforward to see that if $\mu_0(L)$ is not injective then $h^0(L),
h^0(K_C\otimes L^{-1})\geq 2$ and it is an old problem to describe
the locus in $\cM_g$ where the Gieseker-Petri Theorem fails, in particular to
determine its components and their dimensions.
\newline
\indent We fix integers $r, d\geq 1$ such that
$\rho(g,r,d)=g-(r+1)(g-d+r)\geq 0$. As usual, $G^r_d(C)$ is the
variety of linear systems $\mathfrak g^r_d$ on $C$, and if
$(L,V)\in G^r_d(C)$,  we denote by $\mu_0(V):V\otimes
H^0(K_C\otimes L^{-1})\rightarrow H^0(K_C)$ the multiplication
map. We define the Gieseker-Petri locus of type $(r,d)$
$$GP^r_{g,d}:=\{[C]\in \cM_g:\exists \mbox{ a base point free }(L,V) \in G^r_d(C)
\mbox{ with } \mu_0(V) \mbox{  not injective}\}.$$ There
are only two instances when this locus in well understood. First,
$GP^1_{g,g-1}$ can be identified with the above introduced locus
$\cM_g^1$ of curves with a vanishing theta-null which is known to
be an irreducible divisor (cf. \cite{T3}). Then for even $g\geq 4$, $GP^1_{g,(g+2)/2}$ is a
divisor on $\cM_g$ which has an alternate description as the
branch locus of the natural map $H_{g,(g+2)/2}\rightarrow \cM_g$
from the Hurwitz scheme of coverings of $\PP^1$ of degree
$(g+2)/2$ with source curve of genus $g$. This last divisor played
a crucial role in the proof that $\cM_g$ is of general type for
even $g\geq 24$ (cf. \cite{EH3}). It is natural to ask whether
more generally, all loci $GP^r_{g,d}$ are divisors and we give a
partial affirmative answer to this question:
\begin{theorem}\label{petri}
For integers $g\geq 4$ and $(g+2)/2\leq k\leq g-1$, the
Giseker-Petri locus $GP^1_{g,k}$ has a divisorial component.
\end{theorem}

As an easy consequence we mention the following:
\begin{corollary}\label{pol}
For $g\geq 4$ and $0\leq n\leq g-4$, the moduli space $\S_{g,n}^1$
has at least one component of dimension $3g-4$.
\end{corollary}
This last statement can be compared to Polishchuk's recent result
that the moduli space $\S_{g,n}^0$ is of pure dimension $3g-3+n/2$
(cf. \cite{Po}, Theorem 1.1).

\section{Limit theta-characteristics}
In this section, after briefly recalling some basic facts about stable spin curves, we characterize limit theta-characteristics on certain stable curves of compact type after which we prove Theorem \ref{th2}.
\newline
\indent  We review a few things about the moduli space $\ss_g$ (see \cite{C} for more details). If $X$ is a nodal curve, a smooth rational component $R$ of $X$ is called \emph{exceptional} if $\#\bigl(R\cap\overline{(X-R)}\bigr)= 2$. The curve $X$ is called \emph{quasistable} if every two exceptional components are disjoint. Every quasistable curve is obtained by blowing-up some of the nodes of a stable curve.
\newline
\indent A \emph{stable spin curve} consists of a triple $(X,L, \alpha)$, where $X$ is a quasistable curve with $p_a(X)=g$, $L$ is a line bundle on $X$ of degree $g-1$ with $L_R=\OO_R(1)$ for each exceptional component $R$ and $\alpha:L^2\rightarrow \omega_X$ is a homomorphism such that $\alpha_{C}\neq 0$ for any non-exceptional component $C$ of $X$. A \emph{family of stable spin curves}
 is a triple $(f:\cC\rightarrow T, \L,\alpha)$, where $f:\cC\rightarrow T$ is a flat family of quasistable curves, $\L$ is a line bundle on $\cC$ and $\alpha:\L^2\rightarrow \omega_f$ is a homomorphism such that $\alpha_{C_t}$ gives a spin structure on each fibre $C_t=f^{-1}(t)$.
\newline
\indent The stack $\ss_g$ of stable spin curves of genus $g$ has been constructed in \cite{C} where it is also proved that there exists a finite map $\pi:\ss_g\rightarrow \mm_g$ whose fibre over $[C]\in \mm_g$ is the set of stable spin structures on quasistable curves stably equivalent to $C$.
\begin{remark}
Suppose $C=C_1\cup_p C_2$ is a curve of compact type
with $C_1$ and $C_2$ being smooth curves and $g(C_1)=i, g(C_2)=g-i$.
 Then it is easy to see that there are no spin structures on $C$ itself.
 In fact, $\pi^{-1}([C])$ consists of spin structures on the quasistable curve $X=C_1\cup_q R\cup_r C_2$ obtained from $C$ by \lq \lq blowing-up'' $C$ at the node $p$. Each such spin structure is given by a line bundle $L$ on $X$ such that $L_{C_1}^2=K_{C_1}$, $L_{C_2}^2=K_{C_2}$ and $L_R=\OO_R(1)$. More generally, a spin structure on any curve of compact type corresponds to a collection of theta-characteristics on the components.
\end{remark}

%\begin{proposition}\label{ha}
%Let $(f:\cC\rightarrow T, \L,\alpha:\L^2\rightarrow \omega_f)$ be a family of spin curves of compact type. We set $C_t=f^{-1}(t)$ and $L_t=\L_{C_t}$ for $t\in T$. Then
%\begin{enumerate}
%\item the function $T\ni t\mapsto h^0(C_t,L_t) \mbox{ mod }2$ is locally constant.
%\item the locus $T^r=\{t\in T:h^0(C_t,L_t)\geq r+1\}$ is either empty or has codimension $\leq r(r+1)/2$ in $T$.
%\end{enumerate}
%\end{proposition}
Assume now that  $C=C_1\cup_p C_2$ is a curve of compact type
where $C_1$ and $C_2$ are smooth curves of genus $i$ and $g-i$
respectively.
 We define an \emph{r-dimensional limit theta-characteristic} on $C$ (in short, a limit $\theta_g^r$),
 as being a pair of
line bundles $(L_1,L_2)$ with $L_i\in \mbox{Pic}^{g-1}(C_i)$,
together with $(r+1)$-dimensional subspaces  $V_i\subset H^0(L_i)$
such that
\begin{enumerate}
\item $\{l_i=(L_i,V_i)\}_{i=1,2}$ is a limit linear series
$\mathfrak g_{g-1}^r$ in the sense of \cite{EH1}. \item $
L_1^2=K_{C_1}(2(g-i)p)$ and $L_2^2=K_{C_2}(2ip)$.
\end{enumerate}

Using this terminology we  now characterize singular curves in $\mm_{g}^r$:
\begin{lemma} Suppose $[C=C_1\cup_p C_2]\in \mm_g^r$. Then $C$ possesses a $\theta_g^r$.
\end{lemma}
\begin{proof} We may assume that there exists a $1$-dimensional family of
 curves $f:\cC\rightarrow B$ with smooth general fibre $C_b$ and central fibre $C_0=f^{-1}(0)$ stably equivalent to $C$, together with  a line bundle $\L$ on $\cC-C_0$ and a rank $(r+1)$ subvector bundle $V\subset f_*(\L)$ over $B^*:=B-\{0\}$ such that $\L^2_{C_b}\equiv \omega_{C_b}$ for all $b\in B^*$. Then for $i=1,2$ there are unique line bundles $\L_i$ on $\cC$ for extending $\L$ and such that $\mbox{deg}_Y(\L_i)=0$ for every component $Y$ of $C_0$ different from $C_i$. If we denote by $L_i:=\L_i{|_{C_i}}$ and $V_i\subset H^0(L_i)$ the $(r+1)$-dimensional subspace of sections that are limits in $L_i$ of sections in $V$, then by Theorem 2.6 of \cite{EH1} we know that $\{(L_i,V_i)\}_{i=1,2}$ is
a limit $\mathfrak{g}_{g-1}^r$. Finally, since $L_i^2$ and $K_{C_i}$ are isomorphic off $p$ they must differ by a divisor supported at $p$ which accounts for condition $(2)$ in the definition of a $\theta_g^r$.
\end{proof}

We describe explicitly the points in $\mm_g^r\cap \Delta_1$, where $\Delta_1$ is the divisor of curves with an elliptic tail:
\begin{proposition} \label{d1}
Let $[C=C_1\cup _p E]$ be a stable curve with $C_1$ smooth of genus $g-1$ and $E$ an elliptic curve. If $[C]\in \mm_g^r$ then either (1) $[C_1]\in \cM_{g-1}^r$, or (2) there exists a line bundle $L_1$ on $C_1$ such that $(C_1,L_1)\in \S_{g-1}^{r-1}$ and $p\in \rm{Bs}$$|L_1|$. If moreover $p\in C_1$ is a general point, then possibility (2) does not occur hence $[C_1]\in \cM_{g-1}^r$.
\end{proposition}
\begin{proof} We know that $C$ carries a limit $\theta_g^r$, say $l=\{l_{C_1},l_E\}$. By the compatibility
relation between $l_{C_1}$ and $l_E$, the vanishing sequence
$a^{l_{C_1}}(p)$ of $l_{C_1}$ at $p$ is $\geq (0,2,\ldots,r+1)$.
If $l_{C_1}$ has a base point at $p$ then if we set
$L:=L_{C_1}(-p)$ we see that $(C_1,L)\in \S_{g-1}^r$ and we are in
case (1). Otherwise we set $M:=L_{C_1}(-2p)$ and then
$h^0(C_1,M)=r, M^2=K_{C_1}(-2p)$ and $|M+p|$ is a
theta-characteristic on $C_1$ having $p$ as a base point.
\newline
\indent For the  the last statement, we note that a curve has
finitely many positive dimensional theta-characteristics each of
them having only a finite number of base points, so possibility
(2) occurs for at most finitely many points $p\in C_1$.
\end{proof}
We can now  prove Theorem \ref{th2}. More precisely we have the
following result:
\begin{proposition} Fix $r,g\geq 1$. If $\S_{g-1}^r$ has a component of codimension ${r+1 \choose 2}$ in $\S_{g-1}$, then $\S_g^r$ has a component of codimension ${r+1 \choose 2}$ in $\S_g$.
\end{proposition}
\begin{proof}
Suppose $[C_1,L_1]\in \S_{g-1}^r$ is a point for which there
exists a component $\mathcal{Z} \ni [C_1,L_1]$ of $\S_{g-1}^r$
with $\mbox{codim}(\mathcal{Z},\S_{g-1})={r+1 \choose 2}$. We fix a
general point $p\in C_1$ and set $C:=C_1\cup_p E$, where $(E,p)$
is a general elliptic curve. We denote by $X:=C_1\cup_q R\cup_s E$
the curve obtained from $C$ by blowing-up $p$, and we construct a
spin structure on $X$ given by a line bundle $L$ on $X$ with
$L_{C_1}=L_1, L_R=\OO_R(1) $ and $L_E=\OO_E(t-s)$, where $t-s$ is
a non-zero torsion point of order $2$. Clearly
$h^0(X,L)=h^0(C_1,L_1)\geq r+1$. We first claim that $(X,L)$ is a
smoothable spin structure which would show that $[X,L] \in
\ss_g^r$.
\newline
\indent To see this we denote by $(f:\mathcal{X}\rightarrow B, \L, \alpha:\L^2\rightarrow \omega_f)$ the
 versal deformation space of $(X,L)$, so that if $B_1$ denotes the versal deformation space of the
 stable model $C$ of $X$, there is a commutative diagram:

$$\begin{array}{ccccc}

     B &   \stackrel{\sigma}  \longrightarrow & {B}/{\mbox{Aut}(X,L)} & \hookrightarrow  & \ss_g \\
      \rmapdown{\phi} & \; & \; & \;  & \rmapdown{\pi}  \\
    B_1 &    \longrightarrow & {B_1}/{\mbox{Aut}(C)} & \hookrightarrow  & \mm_g \\
\end{array}$$

We define $B^r:=\{b\in B:h^0(X_b,L_b)\geq r+1, h^0(X_b,L_b)\equiv
r+1 \mbox{ mod }2\}$ and Theorem 1.10 from \cite{H} gives that
every component of $B^r$ has dimension $\geq
\mbox{dim}(B)-r(r+1)/2$. We also consider the divisor $\Delta
\subset B$ corresponding to singular spin curves. To conclude that
$(X,L)$ is smoothable we show that there exists a component $W\ni
0$ of $B^r$ not contained in $ \Delta$ (here $0\in B$ is the point
corresponding to $(X,L)$).
\newline
\indent Assume that on the contrary, every component of $B^r$ containing $0$ sits inside $\Delta$. It is straightforward to describe $B^r\cap \Delta$: if $(X_b=C_b\cup R_b\cup E_b,L_b)$ where $b\in B, g(C_b)=g-1,g(E_b)=1$, is a spin curve with
$h^0(X_b,L_b)\geq r+1$, then either (1) $h^0(C_b,L_{b |C_b})\geq r+1$ or (2) $h^0(C_b,L_{b |C_b})=r$ and $L_{b |E_b}=\OO_{E_b}$ (put it differently, $L_{b |E_b}$ is the only odd theta characteristic on $E_b$). Since even and odd theta characteristics do not mix, it follows that any component $0\in W\subset B^r$ will consist entirely of elements $b$ for which $h^0(C_b,L_{b |C_b})\geq r+1$. Moreover, there is a 1:1 correspondence between such components of $B^r$ and components of $\S_{g-1}^r$ through $[C_1,L_1]$. But then the locus
$$\mathcal{Z}_1:=\{b\in \Delta:[C_b,L_{b|C_b}]\in \mathcal{Z}, h^0(E_b,L_{b|E_b})=0\}$$
is a component of $B^r$ containing $0$ and $\mbox{dim}(\mathcal{Z}_1)=\mbox{dim}(\mathcal{Z})+2=3g-4-{r+1 \choose 2},$
 which contradicts the estimate on $\mbox{dim}(B^r)$.

Thus $(X,L)$ is smoothable. We now show that at least one component of $\ss_g^r$ passing through $[C,L_C]$ has codimension ${r+1 \choose 2}$. Suppose this is not the case. Then each component of $\ss_g^r\cap \sigma(\Delta)$ through $[C,L_C]$ has codimension $\leq {r+1 \choose 2} -1$ in $\sigma(\Delta)$. Recalling that $p\in C_1$ was general, Proposition \ref{d1} says that any such component corresponds to curves $C_1'\cap E'$ where $E'$ is elliptic and $[C_1']\in \cM_{g-1}^r$. But then $\sigma(\mathcal{Z}_1)$ is such a component and we have already seen that $\mbox{codim}(\mathcal{Z}_1,\Delta)={r+1 \choose 2}$, which yields the desired contradiction.
\end{proof}

\begin{remark} Retaining the notation from the proof of Theorem \ref{d1}, if
$[C_1,L_1]\in \S_{g-1}^r$ is such that $L_1$ is very ample, then a
smoothing $[C',L_{C'}]\in \S_g^r$ of $[C=C_1\cup_p E, L_C]$
corresponds to a very ample $L_{C'}$. Indeed, assuming by
contradiction that there exist points $x,y\in C'$ such that
$h^0(L_{C'}(-x-y))\geq h^0(L_{C'})-1$, we have three possibilities
depending on the position of the points $r,s\in C$ to which $x$
and $y$ specialize. The case $x,y\in E$ can be ruled out
immediately, while  $x,y\in C_1$ would contradict the assumption
that $L_1$ is a very ample line bundle. Finally, if $x\in C_1$ and
$y\in E$, one obtains that $\{x,p\}$ fails to impose independent
conditions on $|L_1|$, a contradiction. Thus $L_{C'}$ is very
ample.
\end{remark}

\section{Gaussian maps and theta-characteristics}

It may be helpful to review a few things about Gaussian maps on
curves and to explain the connection between Gaussians and
theta-characteristics. This will enable us to construct components
of $\S_g^r$ of dimension achieving the Harris bound.
\newline
\indent For a smooth projective variety $X$ and a line bundle $L$,
we denote by $R(L)$ the kernel of the multiplication map
$H^0(L)\otimes H^0(L)\rightarrow H^0(L^2)$. Following J. Wahl (see
e.g. \cite{W1}), we consider the \emph{Gaussian map}
$\Phi_L=\Phi_{X,L}:R(L)\rightarrow H^0(\Omega_X^1\otimes L^2)$,
defined locally by
$$ s\otimes t \mapsto s\ dt-t\ ds.$$
Since $R(L) =\wedge^2H^0(L)\oplus S_2(L)$,
 with $S_2(L)=\mbox{Ker}\{\mbox{Sym}^2H^0(L)\stackrel{\mu_L} \longrightarrow H^0(L^2)\}$,
  it is clear that $\Phi_L$ vanishes on symmetric tensors  and it makes sense to look at he  restriction
$$\psi_L=\psi_{X,L}:=\Phi_{L | \wedge^2 H^0(L)}:\wedge^2 H^0(L )\rightarrow H^0(\Omega_X^1\otimes L^2).$$

\noindent If $X\subset \PP^r$ is an embedded variety with
$L=\OO_X(1)$, one has the following interpretation for $\Phi_L$:
we  pull back the Euler sequence to $X$ to obtain that
$R(L)=H^0(\Omega_{\PP^r | X}^1 \otimes L^2)$ and then  $\Phi_L$
can be thought of as the map obtained by passing to global
sections in the morphism $\Omega_{\PP^r |X}^1\otimes
L^2\rightarrow \Omega_X^1\otimes L^2.$ Furthermore, if $N_X$ is
the normal bundle of $X$ in $\PP^r$, tensoring the exact sequence
\begin{equation} \label{norm}
0\longrightarrow N_X^{\vee}\longrightarrow \Omega_{\PP^r |X}^1\longrightarrow \Omega_X^1\longrightarrow 0
\end{equation}
by $\OO_X(2)$, we obtain that $\mbox{Ker}(\Phi_L)=\mbox{Ker}(\psi_L)\oplus S_2(L)=H^0(N_X^{\vee}(2))$. If $X$ is projectively normal, from the exact sequence $0\rightarrow \cI_X^2\rightarrow \cI_X\rightarrow N_X^{\vee}\rightarrow 0$
it is straightforward to
check that $\mbox{Ker}(\psi_L)=H^1(\PP^r,\cI_X^2(2)).$

The map $\psi_L$ has been extensively studied especially when $X$
is a curve, in the context of the deformation theory of the cone
over $X$ (cf. e.g. \cite{W1}). The connection between Gaussian
maps and spin curves is given by the following tangent space
computation due to Nagaraj (cf. \cite{N}, Theorem 1): for
$(C,L)\in \S_g^r$, if we make the standard identifications
$T_{[C,L]}(\S_g)=T_{[C]}(\cM_g)=H^1(C, T_C)=H^0(C,K_C^2)^{\vee}$,
then
$$T_{[C,L]}(\S_g^r)=\Bigl(\mbox{Im}(\psi_L):\wedge^2 H^0(L)\rightarrow H^0(K_C^2)\Bigr)^{\perp}.$$
In other words, to show that a component
 $\mathcal{Z}$ of $\S_g^r$ has codimension ${r+1 \choose 2}$ in $\S_g$,
  it suffices to exhibit a spin curve $[C,L]\in \mathcal{Z}$
   such that $h^0(L)=r+1$ and $\psi_L$ is injective. We construct such curves
   as sections of certain homogeneous spaces having injective Gaussians and then we apply
    Theorem \ref{th2} to increase the range of $(g,r)$ for which we have a component of $\S_g^r$ of
    codimension ${r+1 \choose 2}$. We will use repeatedly the following result of  Wahl relating the Gaussian map of a
variety to that of one of its sections (cf. \cite{W2},
Propositions 3.2 and 3.6):
\begin{proposition}\label{wahl}
1. Suppose $X\subset \PP^r$ is a smooth, projectively normal
variety such that $\psi_{X, \OO_X(1)}$ is injective.  If $Y\subset X$ is  a subvariety with ideal sheaf
$\cI$ satisfying the conditions
$$H^1(X, \cI(1))=0,\  H^1(X,\cI^2(2))=0,  \ H^1(X, N_X^{\vee}(2)\otimes \cI)=0,$$
then the Gaussian $\psi_{Y, \OO_Y(1)}$ is injective too.
\newline
\indent 2.  Let $X\subset \PP^r$ be a smooth,  projectively
normal,  arithmetically Cohen-Macaulay variety and $Y=X\cap
\PP^{r-n}\subset \PP^{r-n}$ a smooth codimension $n$ linear
section, where $n<r$. If $H^i\bigl(X,N_X^{\vee}(2-i)\bigr)=0$ for
$1\leq i\leq n$ and $\psi_{X,\OO_X(1)}$ is injective,  then $Y$ is projectively normal and the
Gaussian $\psi_{Y,\OO_Y(1)}$ is also  injective.
\end{proposition}

We will  apply Proposition \ref{wahl} in the case of the
Grassmannian  $X=G(2,n)$ of $2$-dimensional quotients of $\mathbb
C^n$ and for the line bundle $L=\OO_{G(2,n)}(1)$  which gives the
Pl\"{u}cker embedding. In this case $\psi_{\OO_{G(2,n)}(1)}$ is
bijective (cf. \cite{W2}, Theorem 2.11).

 We need to compute the
cohomology of several vector bundles on $G(2,n)$ and we do this
using Bott's theorem (see \cite{FH} for a standard reference). Recall that $G(2,n)=SL_n(\mathbb C)/P$,
where the reductive part of the parabolic subgroup $P$ consists of
matrices of type $\mbox{diag}(A,B)\in SL_n(\mathbb C)$ where $A\in
GL_2(\mathbb C)$ and $B\in GL_{n-2}(\mathbb C)$. We denote by
$\cQ$ the universal rank $2$ quotient bundle defined by the
tautological sequence $$0\longrightarrow \cU\longrightarrow
\OO_{G(2,n)}^{\oplus n}\longrightarrow \cQ\longrightarrow 0.$$

Every irreducible vector bundle over $G(2,n)$ comes from a
representation of the reductive part of $P$. If $e_1, \ldots, e_n$
is an orthonormal basis of $\mathbb R^n$, the positive roots of
$SL_n(\mathbb C)$ are $\{e_i-e_j\}_{i<j}$ and we use the notation
$E(a_1,\ldots ,a_n)$ for the vector bundle corresponding to the
representation with highest weight $a_1 e_1+\cdots +a_n e_n$. We
then have the identifications $\cQ=E(1,0,\ldots ,0)$,
$\OO_{G(2,n)}(1)=\mbox{det}(\cQ)=E(1,1,0,\ldots ,0)$ and
$\cU=E(0,0,1,0,\ldots ,0).$ The cotangent bundle
$\Omega_{G(2,n)}^1=\cQ^{\vee}\otimes \cU$ is irreducible and
corresponds to the highest weight $(0,-1,1,0,\ldots, 0)$. Bott's
theorem can be interpreted as saying that the cohomology group
$H^i\bigl(G(2,n), E(a_1,\ldots ,a_n)\bigr)$ does not vanish if and
only if $i$ is the number of strict inversions in the sequence
$(n+a_1, n-1+a_2, \ldots, 1+a_n)$ and all the entries of this
sequence are distinct.

 First we
establish the following vanishing result:
\begin{proposition}\label{grass}
Let $\GG=G(2,n)\subset \PP^N$ with $ N={n \choose 2}-1$, be the
Grassmannian of lines in its Pl\"{u}cker embedding. We have the
following vanishing statements: \begin{enumerate} \item $H^i(
N_{\GG}^{\vee}(2-i))=0$ for all $ 1\leq i \leq 2n-5, i\neq 2$ and
for $i=2$ and $n\leq 6$. \item $H^i(\Omega_{\GG}^1\otimes
\cQ(-i))=0$ for $0\leq i\leq 2n-7$. \item $H^{i+1}(N_{\GG}^{\vee}
\otimes \cQ(-i)\bigr)=0$ for $1\leq i\leq \rm{min}$$(n,2n-7).$
 \item $H^*(\cQ(-i))=0$ for $1\leq i\leq n$.
\item $H^{i+1}(N_{\GG}^{\vee}(-i))=0$ for $0\leq i\leq n-1$.
\end{enumerate}
\end{proposition}
\begin{proof} (1). We start with the case $i\geq 3$. From the  exact
sequence (\ref{norm})  it suffices to show that (a) $H^{i-1}(\GG,
\Omega_{\GG}^1(2-i))=0$ and that (b) $H^i\bigl(\GG, \Omega_{\PP^N
|\GG}^1(2-i)\bigr)=0$. From the Euler sequence  (b) at its turn is
implied by the vanishings
$H^{i-1}(\OO_{\GG}(2-i))=H^i(\OO_{\GG}(1-i))=0$ which are obvious,
while (a) is a consequence of Bott's theorem (or of Kodaira -
 Nakano vanishing).
 When $i=1$, one checks that $H^0\bigl(\GG,
\Omega_{\GG}^1(1)\bigr)=0$ (Bott again), and that $H^1(\GG,
\Omega_{\PP^N |\GG}^1(1))=0$ (Euler sequence). The remaining case
$i=2$ is handled differently and we employ the Griffiths vanishing
theorem: since $\GG$ is scheme theoretically cut out by quadrics,
the vector bundle $E=N_{\GG}^{\vee}(2)$ is globally generated.
From the exact sequence (\ref{norm}) one finds that
$\mbox{det}(E)=\OO_{\GG}((n-3)(n-4)/2)$ and  we can write
$N_{\GG}^{\vee}=K_{\GG}\otimes E\otimes \mbox{det}(E)\otimes L$,
with $L$ an ample line bundle, precisely when $ n\leq 6$.

Part (2) is a consequence of Le Potier vanishing (cf. \cite{LP}),
while (4) follows from Bott vanishing since
$\cQ(-i)=E(1-i,-i,0,\ldots ,0)$. To prove (3) we tensor the exact
sequence (\ref{norm}) by $\cQ(-i)$ and we have to show that
$H^i(\Omega_{\GG}^1\otimes \cQ(-i))=H^{i+1}(\Omega^1_{\PP^N |\GG}
\otimes \cQ(-i))=0$ which we already treated in parts (2) and (4).
Finally, (5) is handled similarly to (1) and we omit the details.
\end{proof}

For certain $r$ we  construct half-canonical curves $C\subset \PP^r$ of genus $g(r)$ with injective Gaussian. This combined with  Theorem \ref{th2}  proves Theorem \ref{th1}.
\begin{proposition}
For $3\leq r\leq 11, r\neq 10$, there exists a smooth half-canonical curve $C\subset \PP^r$ of genus $g(r)$ (to be specified in the proof), such that the Gaussian map $\psi_{\OO_C(1)}$ is injective. It follows that $\S_{g(r)}^r$ is smooth of codimension $r(r+1)/2$ at the point $[C, \OO_C(1)]$.
\end{proposition}
\begin{proof}
Each case will require a different construction. We treat every
situation separately in increasing order of difficulty.
\newline
\noindent {\textbf{r=3.}}  We let $C$  be a $(3,3)$ complete
intersection in $\PP^3$, hence $g(C)=g(3)=10$ and $K_C=\OO_C(2)$.
Clearly $N_C=\OO_C(3)\oplus \OO_C(3)$, so trivially
$H^1(N_C^{\vee}(2))=0$ which proves that $\psi_{\OO_C(1)}$ is
injective.

\noindent {\textbf{r=4.}} Now $C$ is a complete intersection of
type $(2,2,3)$ in $\PP^4$. Then $g(C)=g(4)=13$ and
$N_C=\OO_C(2)^2\oplus \OO_C(3)$. Using that $C$ is projectively
normal we get that $H^1(\PP^4,\cI_C^2(2))=0$, hence
$\psi_{\OO_C(1)}$ is injective again.

\noindent {\textbf{r=5.}} This is the last case when $C$ can be a
complete intersection: $C$ is of type $(2,2,2,2)$ in $\PP^5$, thus
$g(C)=g(5)=17$ and like in the $r=4$ case we check that
$H^1(\PP^5,\cI_C^2(2))=0$.

\noindent {\textbf{r=8.}} We choose the Grassmannian
$G(2,6)\subset \PP^{14}$. A general codimension $6$ linear section
of $G(2,6)$ is a $K3$ surface $S\subset \PP^8$ with $\mbox{deg}
(S)=14$ and we let $C:=S\cap Q\subset \PP^8$  be a quadric section
of $S$. Then $C$ is half-canonical and $g(C)=g(8)=29$. We claim
that $\psi_{S,\OO_S(1)}$ is injective, which follows from
Proposition \ref{wahl} since
$H^i\bigl(N_{G(2,6)}^{\vee}(2-i)\bigr)=0$ for $1\leq i\leq 6$. To
obtain that $H^1(\PP^8, \cI_C^2(2))=0$, by Proposition \ref{wahl}
we have to check that $H^1(S,\OO_S(-1))=H^1(S,\OO_S(-2))=0$
(Kodaira vanishing), and that $H^1(N_S^{\vee})=0 \Leftrightarrow
H^1(N_S)=0$. Note that  $S$ is a general $K3$ surface of genus $8$
having $\rho(S)=1$ and since by transcendental theory, the Hilbert
scheme of such $K3$ surfaces is irreducible, it will suffice to
exhibit a single $K3$ surface of genus $8$ having this property:
we let $S$ degenerate to a union $R_1\cup_B R_2$ of two rational
scrolls of degree $7$ in $\PP^8$ joined along an elliptic curve
$B\in |-K_{R_i}|$ for $i=1,2$. Then $R_1\cup_B R_2$ is a limit of
smooth $K3$ surfaces $X\subset \PP^8$ of degree $14$ and
$H^1(R_1\cup_B R_2,N_{R_1\cup _B R_2})=0$ (see \cite{CLM}, Theorem
1.2 for more details on this degeneration). It follows that
$H^1(X, N_X)=0$, for a general prime $K3$ surface $X\subset \PP^8$
of genus degree $14$ and then $H^1(S,N_S)=0$ as well.

\noindent {\textbf{r=7.}} In this situation we choose the
$10$-dimensional spinor variety $X\subset \PP^{15}$ corresponding
to a half-spin representation of $\mbox{Spin}(10)$ (see \cite{M}
for a  description of the projective geometry of $X$). One has
that $X$ is a homogeneous space for $SO(10)$, $K_X=\OO_X(-8)$ and
$\mbox{deg}(X)=12$. A general codimension $8$ linear section of
$X$ is a $K3$ surface $S\subset \PP^7$ of degree $12$. Take now
$C$ to be a quadric section of $S$ and then $K_C=\OO_C(2)$ and
$g(C)=g(7)=25$. Since  $N_X^{\vee}$ is irreducible (cf. e.g.
\cite{W2}, Theorem 2.14),  we obtain that the Gaussian map
$\psi_{X, \OO_X(1)}$ is injective.
\newline
\indent To show that $\psi_{S,\OO_S(1)}$ is injective we verify
that $H^i\bigl(N_{X}^{\vee}(2-i)\bigr)=0$ for $1\leq i\leq 8$. For
$3\leq i\leq 8$ this follows from Kodaira-Nakano vanishing for the
twists of sheaves of holomorphic forms on $X$ in a way similar to
the proof of Proposition \ref{grass}, while the $i=1$ it is a
consequence of Bott vanishing. For $i=2$ we use Griffiths
vanishing: since $X$ is cut out by quadrics (see e.g. \cite{M},
Proposition 1.9), the vector bundle $E:=N_X^{\vee}(2)$ is globally
generated, $\mbox{det}(E)=\OO_X(2)$ and one can write
$N_X^{\vee}=K_X\otimes E\otimes \mbox{det}(E)\otimes \OO_X(4)$. In
this way we obtain that $H^2(N_X^{\vee})=0$. Thus $\psi_{S,
\OO_S(1)}$ is injective, and to have the same conclusion for the
Gaussian of $C$, the only non-trivial thing to check is that
$H^1(N_S)=0$, which can be seen by letting $S$ degenerate again to
a union of two rational scrolls like in the case $r=8$.

\noindent {\textbf{r=6.}} We consider the Grassmannian
$\GG=G(2,5)\subset \PP^9$ and we denote by $X\subset \PP^6$ a general codimension
$3$ linear section of $\GG$, by $S:=X\cap Q$ a general quadric section of $X$ and
 by $C:=S\cap Q'$ a general quadric section of $S$. Then $S$ is a $K3$ surface of
 genus $6$, $K_C=\OO_C(2)$ and $g(C)=g(6)=21$. Using Propositions \ref{wahl} and \ref{grass}
 we see easily that $\psi_{X,\OO_X(1)}$ is injective. We claim that $\psi_{S, \OO_S(1)}$ is
 injective as well which would follow from $H^1(X, N_X^{\vee})=0$.
 Since $N_{X/\GG}^{\vee}=\OO_X(-1)^{{\oplus} 3}$, the vanishing of
 $H^1(X, N_X^{\vee})$ is implied by that of $H^1(N_{\GG}^{\vee}\otimes \OO_X)$
 which in its turn is implied by $H^{i+1}(N_{\GG}^{\vee}(-i))=0$ for $0\leq i\leq 3$ (use the Koszul resolution).
 These last vanishing statements  are contained in Proposition \ref{grass} and in this way we obtain
 that $\psi_{S, \OO_S(1)}$ is injective. We finally descend to $C$. To conclude that $\psi_{C, \OO_C(1)}$ is
 injective it is enough to verify that $H^1(N_S)=0$. We could check this  again via the Koszul complex, but
 it is more economical to use that $S$ is a general $K3$ surface of genus $6$ and to invoke
 once more \cite{CLM}, Theorem 1.2, like in the previous cases.

\noindent {\textbf{r=11.}} We start with the Grassmannian
$X=G(2,7)\subset \PP^{20}$ for which $K_X=\OO_X(-7)$ and we let
$C$ be a general codimension $9$ linear section of $X$. Then
$C\subset \PP^{11}$ is a smooth half-canonical curve  of genus
$g(C)=g(11)=43$. To conclude that $\psi_{C,\OO_C(1)}$ is injective
we apply directly the second part of Proposition \ref{wahl}: the
vanishing $H^i(N_{G(2,7)}^{\vee}(2-i))=0$ is guaranteed by
Proposition \ref{grass} for all $1\leq i\leq 9, i\neq 2$. For
$i=2$ we can no longer employ Griffiths vanishing so we proceed
differently: we use (\ref{norm}) together with the vanishing
$H^2(X, \Omega ^1 _{\PP^{20} |X})=0$ coming from the Euler
sequence, to write down the exact sequence
\begin{equation} \label{final}
0\longrightarrow H^1(N_X^{\vee})\longrightarrow H^1(\Omega ^1
_{\PP^{20} |X})\longrightarrow H^1(\Omega_X^1)\longrightarrow
H^2(N_X^{\vee})\longrightarrow 0,
\end{equation}
where $H^1(\Omega_{\PP^{20} |X}^1)\cong H^0(\OO_X)\cong \mathbb C$
and $H^1(\Omega_X^1)\cong \mathbb C$. From Bott's theorem at most
one of the cohomology groups of the irreducible bundle
$N_X^{\vee}$ are $\neq 0$, hence either $H^2(N_X^{\vee})=0$ and
then we are done, or else, if $H^2(N_X^{\vee})\neq 0$ then
$H^1(N_X^{\vee})=0$, and the map in the middle of the sequence
(\ref{final}) is bijective which yields a contradiction.

\noindent {\textbf{r=9.}} This is the most involved case. We look
at the ample vector bundle $\mathcal{F}:=\cQ(1)$ on
$\GG=G(2,6)\subset \PP^{14}$ and choose a general section $s\in
H^0(\GG, \mathcal{F})$. We denote by $Z$ the zero locus of $s$, by
$\cI=\cI_{Z/\GG}$ the ideal of $Z$ inside $\GG$, and by $\cI_Z$
and $\cI_{\GG}$ the ideals of $Z$ and $\GG$ in $\PP^{14}$
respectively. By adjunction, we have that
$\cI/{\cI^2}=\cQ^{\vee}(-1)\otimes \OO_Z$ and the Koszul complex
gives a resolution for $Z$:
$$0\longrightarrow \OO_{\GG}(-3)\longrightarrow
\cQ^{\vee}(-1)\longrightarrow \cI\longrightarrow 0.$$  We first
claim that $Z\subset \PP^{14}$ is nondegenerate and projectively
normal. This will follow if we show that $H^0(\GG, \cI(1))=0$ and
$H^1(\GG, \cI(r))=0$ for $r\geq 1$. Using the Koszul resolution,
the first vanishing is implied by
$H^0(\cQ^{\vee})=H^1(\OO_{\GG}(-2))=0$ which is clear. For the
second vanishing we have to check that
$H^1(\cQ^{\vee}(r-1))=H^2(\OO_{\GG}(r-3))=0$ for $r\geq 1$. Since
$\cQ^{\vee}(r-1)=E(r-1,r-2,0,0,0,0)$ and
$\OO_{\GG}(r-3)=E(r-3,r-3,0,0,0,0)$ this can be checked instantly
using Bott's theorem. \newline \indent Next we claim that the
$\psi_{Z,\OO_Z(1)}$ is injective. By Proposition \ref{wahl}, we
have to verify that (1) $H^1(\GG, \cI^2(2))=0$ and that (2)
$H^1(\GG, N_{\GG}^{\vee}(2)\otimes \cI)=0$. We start with (1).
From the exact sequence $$0\longrightarrow \cI^2(2)\longrightarrow
\cI(2)\longrightarrow \cQ^{\vee}(1)\otimes \OO_Z \longrightarrow
0,
$$
using that $Z$ is projectively normal, (1) is implied by the
bijectivity of the map $H^0(\cI(2))\rightarrow
H^0(\cQ^{\vee}(1)\otimes \OO_Z)$. This is a consequence of the
isomorphism $\cQ^{\vee}(1)\cong \cQ$ and of the Koszul resolution
giving that $H^0(Z, \cQ^{\vee}(1)\otimes \OO_Z)=H^0(\GG,
\cQ^{\vee}(1))=H^0(\GG, \cI(2))$, where for the first isomorphism
one uses that $H^0(\GG, \cI\otimes \cQ)=H^1(\GG, \cI\otimes
\cQ)=0$, which is straightforward to check via Bott's theorem.
\newline
\indent We turn to (2). The cohomology of $\cI \otimes
N_{\GG}^{\vee}(2)$ is computed from the Koszul complex of $\cI$,
which yields an isomorphism $H^1(N_{\GG}^{\vee}\otimes
\cI(2))=H^1(N_{\GG}^{\vee}\otimes \cQ^{\vee}(1))$ (because we have
$H^i(N_{\GG}^{\vee}(-1))=0$ for $i=1,2$- this being checked via
the sequence (\ref{norm})). Next we write the cohomology sequence
associated to the exact sequence
$$0\longrightarrow N_{\GG}^{\vee}\otimes \cQ^{\vee}(1)\longrightarrow
\Omega^1_{\PP^{14} |\GG}\otimes \cQ^{\vee}(1) \longrightarrow
\Omega_{\GG}^1\otimes \cQ^{\vee}(1)\longrightarrow 0.$$ The map
$H^1(\Omega^1_{\PP^{14} |\GG}\otimes \cQ^{\vee}(1))\rightarrow
H^1(\Omega^1_{\GG}\otimes \cQ^{\vee}(1))$ is an isomorphism: from
the Euler sequence one obtains that $H^1(\Omega^1_{\PP^{14} | \GG}
\otimes \cQ^{\vee}(1))=H^0(\cQ^{\vee}(1))$, while tensoring by
$\Omega_{\GG}^1(1)$ the dual of the tautological sequence, one
gets that $H^1(\Omega^1_{\GG}(1)\otimes
\cQ^{\vee})=H^0(\cU^{\vee}\otimes
\Omega_{\GG}^1(1))=H^0(\cQ^{\vee}(1))$ (or alternatively, use for
this \cite{LP}, Corrolaire 2). Moreover $H^0(\Omega
^1_{\GG}\otimes \cQ^{\vee}(1))$ injects into
$H^0(\Omega^1_{\GG}(1)) ^{\oplus 6}$ which is zero by Bott's
theorem. Hence $H^1(N_{\GG}^{\vee}\otimes \cQ^{\vee}(1))=0$ and
this proves that $\psi_{Z,\OO_{Z}(1)}$ is injective.

We now take a general codimension $5$ linear section of $Z$ which
is a curve $C\subset \PP^9$ with $K_C=\OO_C(2)$. A routine
calculation gives that $\mbox{deg}(C)=3\mbox{deg}(\GG)=42$, hence
$g(C)=g(9)=43$. We claim that $\psi_{C, \OO_C(1)}$ is injective.
Since $\psi_{Z, \OO_{Z}(1)}$ is injective, by Proposition
\ref{wahl} we are left with checking that $Z$ is ACM (this amounts
to $H^i(\OO_Z(j))=0$ for $i\neq 0,6=\mbox{dim}(Z)$, which easily
follows from the Koszul complex) and that
$H^i(Z,N_Z^{\vee}(2-i))=0$ for $1\leq i\leq 5$ (here
$N_Z=(\cI_Z/{\cI_Z})^{\vee}$ is the normal bundle of $Z$ in
$\PP^{14}$). We employ the exact sequence
$$0\longrightarrow N_{\GG}^{\vee}\otimes \OO_Z \longrightarrow N_Z^{\vee} \longrightarrow
\cI/{\cI^2}\longrightarrow 0,$$ from which it will suffice to show
that (a)
$H^i\bigl(Z,\cI/{\cI^2}(2-i)\bigr)=H^i\bigl(\cQ^{\vee}(1-i)\otimes
\OO_Z\bigr)=0$ for $1\leq i\leq 5$ and that  (b)
$H^i(N_{\GG}^{\vee}(2-i)\otimes \OO_Z)=0$, which in turn is a
consequence of
$H^i(N_{\GG}^{\vee}(2-i))=H^{i+1}(N_{\GG}^{\vee}(-1-i))=0$ and of
the vanishing $H^{i+1}(N_{\GG}^{\vee}\otimes \cQ^{\vee}(1-i))=0$
(for all these use Proposition \ref{grass}).
\newline \indent We are left with (a) which is a consequence
of $H^i(\cQ^{\vee}(1-i))=0$ (again, use Proposition \ref{grass}),
of $H^{i+2}\bigl(\cQ^{\vee}(-2-i)\bigr)=0$, and of
$H^{i+1}\bigl(\cQ\otimes \cQ(-2-i)\bigr)=0$. For this last
statement use that $\cQ\otimes \cQ=S^2\cQ\oplus \mbox{det}(\cQ)$
and each summand being irreducible  the vanishing can be easily
verified via Bott's theorem.
\end{proof}

We believe that there should be a uniform way of constructing
half-canonical curves $C\subset \PP^r$ for any $r\geq 3$ of high
genus $g>>r$ and having injective Gaussian maps (though no longer
as sections of homogeneous varieties). Together with Theorem
\ref{th2} this prompts us to make the following:
\begin{conjecture}\label{white}
For any $r\geq 3$ and $g\geq {r+2\choose 2}$, there exists a
component of $\S_g^r$ of codimension ${r+1 \choose 2}$ inside $\S_g$.
\end{conjecture}
The bound $g\geq {r+2\choose 2}$ is obtained by comparing the
expected dimension $3g-3-{r+1\choose 2}$ of $\S_g^r$ with the
expected dimension of the Hilbert scheme $\mbox{Hilb}_{g-1,g,r}$
of curves $C\subset \PP^r$ of genus $g$ and degree $g-1$. We
believe that there exists a component of $\mbox{Hilb}_{g-1,g,r}$
consisting entirely of half-canonically embedded curves. To prove
the Conjecture it would suffice to construct a smooth
half-canonical curve $C\subset \PP^r$ of genus $g={r+2\choose 2}$
such that $H^1(C,N_{C/\PP^r})=0$, that is, $\mbox{Hilb}_{g-1,g,r}$
is smooth at the point $[C]$ and has expected dimension $h^0(C,
N_{C/\PP^r})=4(g-1)$. Note that for such  $C$, the map  $\Psi_{C,
\OO_C(1)}$ would be injective, in particular $C$ would not sit on
any quadrics. This gives the necessary inequality $g\geq {r+2
\choose 2}$.  The main difficulty in proving Conjecture
\ref{white} lies in the fact that the degeneration techniques one
normally uses to construct \lq\lq regular" components of Hilbert
schemes of curves, seem to be at odds with the requirement that
$C$ be half-canonical.

\section{Gieseker-Petri loci}
In this section we construct divisorial components of  the loci $GP^1_{g,k}$.
 The method we use is inductive and close in spirit to the one employed in Section 2 to construct
components of $\S_g^r$ of expected dimension. We begin by
describing a setup that enables us to analyze the following
situation: if $\{L_b\}_{b\in B^*}$ and $\{M_b\}_{b\in B^*}$ are
two families of line bundles over a $1$-dimensional family of
smooth curves $\{X_b\}_{b\in B^*}$, where $B^*=B-\{b_0\}$ with
$b_0\in B$, we want to describe what happens to the multiplication
map
$$\mu_b=\mu_b(L_b,M_b):H^0(X_b, L_b)\otimes H^0(X_b,
M_b)\rightarrow H^0(X_b, L_b\otimes M_b)$$ as $X_b$ degenerates to
a singular curve of compact type $X_0$.

Suppose first that $C$ is a smooth curve and $p\in C$. We recall
that if $l=(L,V)$ is a linear series of type $\mathfrak g^r_d$
with $L\in \mbox{Pic}^d(C)$ and $V\subset H^0(L)$, the
\emph{vanishing sequence} of $l$ at $p$ $$a^l(p):0\leq
a_0^l(p)<\ldots <a^l_r(p)\leq d,$$ is obtained by ordering the set
$\{\mbox{ord}_p(\sigma)\}_{\sigma \in V}$. If  $L$ and $M$  are
line bundles on $C$ and  $\rho \in H^0(L)\otimes H^0(M)$ we write
that $\mbox{ord}_p(\rho)\geq k$, if $\rho$ lies in the span of
elements of the form $\sigma\otimes \tau$, where $\sigma \in
H^0(L)$ and $\tau \in H^0(M)$ are such that
$\mbox{ord}_p(\sigma)+\mbox{ord}_p(\tau)\geq k$.
\newline Let $\mu_{L,M}:H^0(L)\otimes H^0(M)\rightarrow H^0(L\otimes M)
$ be the multiplication map. We shall use the following
observation: suppose $\{\sigma_i\}\subset H^0(L)$ and
$\{\tau_j\}\subset H^0(M)$ are bases of global sections
\emph{adapted to the point} $p\in C$ in the sense that
$\mbox{ord}_p(\sigma_i)=a_i^L(p)$ and
$\mbox{ord}_p(\tau_j)=a_j^M(p)$ for all $i$ and $j$. Then if $\rho
\in \mbox{Ker}(\mu_{L,M})$ then there must exist distinct pairs of
integers $(i_1,j_1)\neq (i_2,j_2)$ such that
$$\mbox{ord}_p(\rho)=\mbox{ord}_p(\sigma_{i_1})+\mbox{ord}_p(\tau_{j_1})=\mbox{ord}_p(\sigma_{i_2})+
\mbox{ord}_p(\tau_{j_2}).$$

Suppose now that $\pi:X\rightarrow B$ is a family of genus $g$ curves over $B={\rm Spec}(R)$,
with $R$ being a complete DVR with local parameter $t$, and let $0, \eta$ denote the special
and the generic point of $B$ respectively. Assume furthermore that $X_{\eta}$ is smooth
 and that $X_0$ is singular but of compact type. If $L_{\eta}$ is a line bundle on $X_{\eta}$
 then, as explained in \cite{EH1}, there is a canonical way to associate to each component $Y$
  of $X_0$ a line bundle $L^Y$ on $X$ such that $\mbox{deg}_Z(L^Y_{|_Z})=0$ for every component
$Z$ of $X_0$ different from $Y$. We set $L_Y:=L^Y_{|_Y}$ which is
a line bundle on the smooth curve $Y$.

We fix $\sigma \in \pi_*L_{\eta}$ a section on the generic fibre. We denote by $\alpha$ the smallest integer such that $t^{\alpha}\sigma \in \pi_*L^Y$, that is, $t^{\alpha}\sigma \in \pi_*L^Y-t\pi_*L^Y$. Then we set
$$\sigma^Y:=t^{\alpha} \sigma \in \pi_*L^Y \mbox{ and } \sigma_Y:=\sigma^Y_{|_Y}\in H^0(Y,L_Y).$$
For a different component $Z$ of the special fibre $X_0$ meeting $Y$ at a point $p$, we define similarly $L^Z,L_Z,\sigma^Z$ and $\sigma_Z$. If  we write $\sigma^Z=t^{\beta}\sigma^Y \in \pi_*L^Z$ for a unique integer $\beta$, we have the following compatibility relation between $\sigma_Y$ and $\sigma_Z$ (cf. \cite{EH1}, Proposition 2.2):
\begin{equation}\label{limit}
\mbox{deg}(L_Y)-\mbox{ord}_p(\sigma_Y)\leq \beta \leq \mbox{ord}_p(\sigma_Z).
\end{equation}
An immediate consequence of this is the inequality $$\mbox{ord}_p(\sigma_Y)+\mbox{ord}_p(\sigma_Z)\geq \mbox{deg}(L_Y)=\mbox{deg}(L_Z).$$

Assume from now on that we have two line bundles $L_{\eta}$ and $M_{\eta}$ on $X_{\eta}$ and we choose an element $\rho \in H^0(X_\eta,L_\eta)\otimes _{R_\eta } H^0(X_\eta, M_\eta).$ If $Y$ and $Z$ are components of $X_0$ meeting at $p$ as above, we define $\rho^Y:=t^{\gamma}\rho \in H^0(X, L^Y)\otimes _R H^0(X,M^Y)$, where $\gamma$ is the minimal integer with this property. We have a similar definition for $\rho^Z\in H^0(X,L^Z)\otimes _R H^0(X,M^Z)$.  Between the sections $\rho^Y$ and $\rho^Z$ there is a relation
$\rho^Z=t^{\alpha} \rho^Y$
 for a uniquely determined integer $\alpha$. To determine $\alpha$ we proceed as follows:
  we choose bases of sections $\{\sigma_i=\sigma^Y_i\}$ for $H^0(X,L^Y)$ and $\{\tau_j=\tau^Y_j\}$ for $H^0(X,M^Y)$
  such that
$\mbox{ord}_p(\sigma_{i, Y})=a^{L_Y}_i(p) \mbox{ and } \mbox{
ord}_p(\tau_{j, Y})=a^{M_Y}_j(p),$ for all relevant $i$ and $j$
(cf. e.g. \cite{EH1}, Lemma 2.3, for the fact that this can be
done). Then there are integers $\alpha_i$ and $\beta_j$ defined by
$\sigma_i^Z=t^{\alpha_i} \sigma_i \mbox{ and }
\tau_j^Z=t^{\beta_j}\tau_j.$ To obtain a formula for the integer
$\alpha$ we write $\rho^Y=\sum_{i,j} f_{ij}\sigma_i \otimes
\tau_j$, where $f_{ij}\in R$. We have the identity
$$\rho^Z=\sum_{i,j} (t^{\alpha-\alpha_i-\beta_j} f_{ij})(t^{\alpha_i} \sigma_i)\otimes (t^{\beta_j} \tau_j),$$
from which we easily deduce that $\alpha=\mbox{max}_{i,j}
\{\alpha_i+\beta_j-\nu(f_{ij})\},$ where $\nu$ denotes the
valuation on $R$ (see also \cite{EH2}, Lemma 3.2).

\begin{lemma}\label{sym}
With the above notations, if $\rho_Y:=\rho^Y_{|_Y}\in H^0(Y,L_Y)\otimes H^0(Y,M_Y)$ and $\rho_Z:=\rho^Z_{|_Z}\in H^0(Z,L_Z)\otimes H^0(Z,M_Z)$, then
\begin{center}
$\rm{ord}$$_p(\rho_Y)+\rm{ord}$$_p(\rho_Z)\geq \rm{deg}$$(L_Y)+\rm{deg}$$(M_Y).$
\end{center}
\end{lemma}
\begin{proof} By definition, there exists a pair of indices $(i_1,j_1)$ such that $\nu(f_{i_1j_1})=0$ and
$$\mbox{ord}_p(\rho_Y)=\mbox{ord}_p(\sigma_{i_1,Y})+\mbox{ord}_p(\sigma_{j_1,Y})$$
and clearly $\alpha\geq \alpha_{i_1}+\beta_{j_1}$. To get an estimate on $\mbox{ord}_p(\rho_Z)$ we only have to take into account the pairs of indices $(i,j)$ for which $\alpha_i+\beta_j=\alpha+\nu(f_{ij})\geq \alpha_{i_1}+\beta_{j_1}$. For at least one such pair $(i,j)$ we have that
$$\mbox{ord}_p(\rho_Z)=\mbox{ord}_p(t^{\alpha_i}\sigma_{i, Z})+\mbox{ord}_p(t^{\beta_j}\tau_{j, Z})\geq \alpha_i+\beta_j.$$

\noindent On the other hand, by applying (\ref{limit}) we can
write
$$\mbox{ord}_p(\rho_Y)=\mbox{ord}_p(\sigma_{i_{1},Y})+\mbox{ord}_p(\tau_{j_1,Y})\geq \mbox{deg}(L_Y)+\mbox{deg}(M_Y)-\alpha_{i_1}-\beta_{j_1},$$
whence we finally have that
$\mbox{ord}_p(\rho_Z)+\mbox{ord}_p(\rho_Y)\geq \mbox{deg}(L_Y)+\mbox{deg}(M_Y).$
\end{proof}

We now fix integers $g$ and $k$ such that $g\geq 4$ and
$(g+2)/2\leq k\leq g-1$ and consider the locus $GP^1_{g,k}$ of
curves $[C]\in \cM_g$ for which the Gieseker-Petri Theorem fails
for a base point free pencil $\mathfrak g^1_k$. We denote by
$\overline{GP}_{g,k}^1$ the closure of $GP^1_{g,k}$ in $\mm_g$ and
we  study $GP^1_{g,k}$ inductively by understanding the
intersection $\overline{GP}^1_{g,k}\cap \Delta_1$.

\begin{definition}
For a smooth curve $C$ of genus $g$, a Gieseker-Petri
\emph{$(gp)^1_{k}$-relation} consists of a linear series $(L,V)\in
G^1_k(C), V\subset H^0(L)$, together with an element
$$\rho \in \PP\mbox{Ker}\{\mu_0(V):V\otimes H^0(K_C\otimes
L^{-1})\rightarrow H^0(K_C)\}.$$ If $C=C_1\cup_p C_2$ is of
compact type with $C_1$ and $C_2$ smooth of genus $i$ and $g-i$
respectively, a \emph{$(gp)_{k}^1$-relation} on $C$ is a
collection $(l,m,\rho_1, \rho_2)$, where $l=\{(L_{C_1},V_{C_1}),
(L_{C_2}, V_{C_2})\}$ is a limit $\mathfrak g^1_k$ on $C$,
$m=\{\bigl(M_{C_1}=K_{C_1}(2(g-i)p)\otimes L_{{C_1}}^{-1},
W_1\bigr), \bigl(M_{C_2}=K_{C_2}(2ip)\otimes L_{C_2}^{-1},
W_2\bigr)\}$ is a limit $\mathfrak g_{2g-2-k}^{g-k}$ on $C$, and
elements $$\rho_1\in \PP\mbox{Ker}\{V_{C_1}\otimes
W_{C_1}\rightarrow H^0\bigl(K_{C_1}(2(g-i)p)\bigr)\}, \rho_{2}\in
\PP \mbox{Ker}\{V_{C_2}\otimes W_{C_2}\rightarrow H^0\bigl(
K_{C_2}(2ip)\bigr)\}$$ satisfying the relation
$\mbox{ord}_p(\rho_1)+\mbox{ord}_p(\rho_2)\geq 2g-2$.
\end{definition}

For a curve $C$ of compact type, we  denote by $Q_{k}^1(C)$ the variety of $(gp)_{k}^1$-relations on  $C$ together with the scheme structure coming from its
natural description as a determinantal variety. The discussion above shows that if $[C]\in \overline{GP}_{g,k}^1$ then $Q_k^1(C)\neq \emptyset$.
 Our strategy is to construct $(gp)_{k}^1$-relations on certain singular curves and prove that they can be deformed to nearby smooth curves filling up a divisor in $\mm_g$. The most important technical result of this section is the construction of
the moduli space of $(gp)_{k}^1$-relations over the versal
deformation space of a curve of compact type inside the divisor
$\Delta_1$:

\begin{theorem}\label{rel}
Fix integers $g\geq 4$ and $k$ such that $(g+2)/2\leq k\leq g-1$.
Let $C$ be a smooth curve of genus $g-1$, $p\in C$ and $X_0:=C\cup
_p E$, where $E$ is an elliptic curve. We denote by
$\pi:X\rightarrow B$ the versal deformation space of $X_0$, with
$X_0=\pi^{-1}(0)$ and $0\in B$. Then there exists a scheme
$\cQ_{k}^1\rightarrow B$, quasi-projective over $B$ and compatible
with base change, such that the fibre over $b\in B$ parametrizes
$(gp)_{k}^1$-relations over $X_b$. Furthermore each component of
$\cQ_{k}^1$ has dimension $\geq \rm{dim}$$(B)-1=3g-4$.
\end{theorem}

\begin{proof} The scheme $\cQ_{k}^1$ is going to be the disjoint
union of subschemes where the vanishing sequences of the aspects
of the two underlying limit linear series of a
$(gp)_{k}^1$-relation are also specified. We will prove the
existence for the component corresponding to vanishing sequences
$(1,2)$ and $(k-2,k-1)$ for the limit $\mathfrak g^1_k$ and
$(1,2,\ldots,g-k+1)$ and $(g-3,g-2,\ldots,2g-3-k)$ for the limit
$\mathfrak g^{g-k}_{2g-2-k}$ respectively. The construction is
entirely similar for the other compatible vanishing sequences. In our
proof we will use  Theorem 3.3 in \cite{EH1} where
a moduli space of limit linear series over the versal deformation
space of a curve of compact type is constructed.
\newline \indent
We start by setting some notations. We denote by $\Delta\subset B$
the \lq \lq boundary" divisor corresponding to curves in which the
node $p$ is not smoothed. We denote by $\cC_p$ and $\cE_p$ the
closures in $X$ of the components of $\pi^{-1}(\Delta)$ containing
$C-\{p\}$ and $E-\{p\}$ respectively. By shrinking $B$ if
necessary we can assume that $\OO_{X}(\cC_p+\cE_p)=\OO_{X}$. We
denote by $\pi_{P^C}:P^C
 \rightarrow B$
 the relative Picard variety corresponding to the family
 $X\rightarrow B$ such that for $b\in \Delta$ and $\pi^{-1}(b)=
X_b=C_b\cup
 E_b$ with $C_b\subset \cC_p$ and $E_b\subset \cE_p$, the fibre of
 $P^C$ over $b$ consists of line bundles $L_b$ on $X_b$ with
 $\mbox{deg}_{C_b}(L_b)=k$ and $\mbox{deg}_{E_b}(L_b)=0$.
 Interchanging the role of $C$ and $E$ we get another Picard
 variety $P^E\rightarrow B$ and tensoring with $\OO_{X}(k\cC_p)$
 gives an isomorphism $P^C\rightarrow P^E$. We denote by $P$ the
 inverse limit of $P^C$ and $P^E$ under this isomorphism. For $b\in
B$ and any line bundle $L$ on $X_b$, we define two new line
bundles $L_C$ and $L_E$ as follows: if $b\in B-\Delta$ then
$L_C=L_E=L$. If $b\in \Delta$ and $X_b=C_b\cup_q E_b$, then
$L_C$ is the restriction to $C$ of the unique line bundle on $X_b$
obtained from $L$ by tensoring with a divisor based at $q$ and
whose restriction to $E_b$ is of degree $0$ (and a similar
definition for $L_E$ with $C$ and $E$ reversed). Proceeding in a a
way identical to \cite{EH1}, pp. 356-360, we  construct a space
of compatible frames $\phi:\F\rightarrow B$ factoring through
$\pi_P:P\rightarrow B$, and which parametrizes objects
 $$ x=\{b, L,(\sigma_i^C)_{i=0,1}, (\sigma_i^E)_{i=0,1},
(\tau_j^C)_{j=0,\ldots, g-k}, (\tau_j^E)_{j=0, \ldots, g-k}\}, $$
where $b\in B$, $L$ is a line bundle of degree $k$ on $X_b$,
$(\sigma_i^C)$ (resp. $(\sigma_i^E)$) is a projective frame inside
$H^0(L_C)$ (resp. $H^0(L_E)$), while $(\tau_j^C)$ (resp.
$(\tau_j^E)$) is a projective frame inside
$H^0\bigl((\omega_{X_b}\otimes L^{-1})_C\bigr)$ (resp.
$H^0\bigl((\omega_{X_b}\otimes L^{-1})_E\bigr))$, subject to the
following identifications: if $b\in B-\Delta$, so $X_b$ is smooth
and $L_E=L_C=L$, then we identify $\sigma_i^C=\sigma_{1-i}^E$ for
$i=0,1$ and $\tau_j^C=\tau_{g-k-j}^E$ for $j=0,\ldots, g-k$ (that
is, there are only two frames, one inside $H^0(L)$, the other
inside $H^0(K_{X_b}\otimes L^{-1})$). If $b\in \Delta$ and
$X_b=C_b\cup_q E_b$ then we require that
$\mbox{ord}_q(\sigma_i^C)\geq i+1$, $\mbox{ord}_q(\sigma_i^E)\geq
k-2+i$ for $i=0,1$, while $\mbox{ord}_q(\tau_j^C)\geq j+1$ and
$\mbox{ord}_q(\tau_j^E)\geq g+j-3$. In this latter case
$l=\{\bigl(L_C,<\sigma_i^C>_i\bigr), \bigl(L_E,
<\sigma_i^E>_i\bigr)\}$ is a limit $\mathfrak g^1_k$ and
$m=\{\bigl((\omega_X\otimes L^{-1})_C, <\tau_j^C>_j\bigr),
\bigl((\omega_X\otimes L^{-1})_E, <\tau_j^E>_j\bigr)\}$ is a limit
$\mathfrak g_{2g-2-k}^{g-k}$ on $X_b$.

The scheme $\F$ is determinantal and each of its components has
dimension $\geq \mbox{dim}(B)+g+2+(g-k+1)(g-k-2)$, which is
consistent with the naive dimension count for the fibre over $b\in
B-\Delta$. We also have tautological line bundles
$\tilde{\sigma}_i^C, \tilde{\sigma}_i^E, \tilde{\sigma}_j^C$ and
$\tilde{\sigma}_j^E$ over $\F$, with fibres over each point being
the $1$-dimensional vector space corresponding to the frame
denoted by the same symbol. For $2\leq i\leq g-k+2$, we consider
the rang $g$ vector bundle bundle
$\Psi_i:=\pi_*\bigl(\omega_{X/B}\otimes \OO_X(i \cC_p)\bigr)$;
hence $\Psi_i(b)=H^0(X_b,L_b)$ for $b\in B-\Delta$, while for
$b\in \Delta$ the fibre $\Psi_i(b)$ consists of those sections in
$H^0\bigl(K_{C_b}(-(i-1) q)\bigr)\oplus H^0\bigl(\OO_{E_b}((i+1)
q)\bigr)$ that are compatible at the node $q$.
\newpage
\indent For $1\leq i\leq g-k$ we define a subscheme $\G_i$ of $\F$
by the equations
\begin{equation}\label{fr}
\tilde{\sigma}_0^C\cdot \tilde{\tau}_i^C=\tilde{\sigma}_1^C\cdot
\tilde{\tau}_{i-1}^C \mbox{ and }\tilde{\sigma}_1^E\cdot
\tilde{\tau}_{g-k-i}^E=\tilde{\sigma}_{0}^E\cdot
\tilde{\tau}_{g-k-i+1}^E.
\end{equation}
Here by $(\tilde{\sigma}_{\alpha}^C\cdot
\tilde{\tau_{\beta}}^C)(x) $ we denote the element in $\PP
H^0((\omega_{X_b})_C)$ obtained by multiplying representatives of
$\tilde{\sigma}_{\alpha}^C(x)$ and of $\tilde{\tau}_{\beta}^C(x)$
for each $x\in \F, b=\phi(x)$. To make more sense of (\ref{fr}),
for each $x\in \F$ the element $\bigl((\tilde{\sigma}_0^C\cdot
\tilde{\tau}_i^C)(x), (\tilde{\sigma}_1^E\cdot
\tilde{\tau}_{g-k-i}^E)(x)\bigr)$ gives rise canonically to a
point in $\PP\bigl((\phi^*\Psi_{i+1})(x)\bigr)$ and abusing the
notation we can consider $(\tilde{\sigma}_0^C\cdot
\tilde{\tau}_i^C, \tilde{\sigma}_1^E\cdot \tilde{\tau}_{g-k-i}^E)$
and $(\tilde{\sigma}_1^C \cdot \tilde{\tau}_{i-1}^C,
\tilde{\sigma}_0^E\cdot \tilde{\tau}_{g-k-i+1}^E)$ as  sections of
the $\PP^{g-1}$ bundle $\PP(\phi^*\Psi_{i+1}))\rightarrow \F$.
Then $\G_i$ is the locus in $\F$ where these sections coincide and
therefore each component of $\G_i$ has dimension $\geq \mbox{dim}(\F)-g+1$.

We define $\cQ_{k}^1$ as the union of the scheme theoretic images
of  $\G_i$ for $1\leq i\leq g-k$ under the map
$$\G_i\ni x\stackrel{\chi_i}\mapsto \bigl(b,l,m,
\rho_1=(\sigma_0^C\otimes \tau_i^C-\sigma_1^C\otimes \tau_{i-1}^C),
\rho_2=(\sigma_1^E\otimes \tau_{g-k-i}^E-\sigma_0^E\otimes
\tau_{g-k-i+1}^E)\bigr),$$ where we recall that $l$ and $m$ denote
the underlying limit $\mathfrak g^1_k$ and $\mathfrak
g_{2g-2-k}^{g-k}$ respectively. From the base point free pencil
trick applied on both $C_b$ and $E_b$, it is easy to see that
$\cQ_{k}^1$ contains all $(gp)_{k}^1$-relations on the curves
$X_b=C_b\cup _q E_b$, the points coming from $\G_i$ corresponding
to those $(b,l,m, \rho_1,\rho_2)$ for which
$\mbox{ord}_q(\rho_1)\geq i+2 $ and $\mbox{ord}_q(\rho_2)\geq
2g-i-4$.

We are left with estimating $\mbox{dim}(\cQ_{k}^1)$: having fixed
$(b,l,m,\rho_1,\rho_2)$ inside $\chi_i(\G_i)$, there are two
cases to consider depending on whether $X_b$ is smooth or not. In
each case we obtain the same estimate for the fibre dimension of
 $\chi_i$ but here we only present the case $b\in \Delta$,
when $X_b=C_b\cup_q E_b$. We have a one dimensional family of
choices for each of $(\sigma_0^C, \sigma_1^C)$ and $(\sigma_0^E,
\sigma_1^E)$, and after choosing these, $(\tau_i^C, \tau_{i-1}^C)$
and $(\tau_{g-k-i}^E, \tau_{g-k-i+1}^E)$ are uniquely determined
(again, use the base point free pencil trick). For choosing the
remaining $\tau_\alpha ^C, \alpha\neq i, i-1$ we have a
$(g-k)(g-k+1)/2-(2g-2k-2i+1)$-dimensional family of possibilities,
while for $\tau_\beta ^E, \beta\neq g-k-i,g-k-i+1$ we get another
$(g-k)(g-k+1)/2-2i+1$ dimensions. Adding these together we get
that each component of $\cQ_{k}^1$ has dimension $\geq 3g-4$.
\end{proof}

We can now prove Theorem \ref{petri}. More precisely we have the
following inductive result:

\begin{theorem}\label{gp}
Fix integers $g,k$ such that $g\geq 4$ and $(g+2)/2\leq k\leq
g-1$. Suppose $GP_{g-1, k-1}^1$ has a divisorial component $Z$ for
which a general  $[C]\in Z$ is such that there exists a
$0$-dimensional component of $Q_{k-1}^1(C)$ whose general point
corresponds to a base point free $\mathfrak g^1_{k-1}$. Then
$GP_{g,k}^1$ has a divisorial component $Z'$, for which a general
curve $[C']\in Z'$ is such that $Q_k^1(C')$ has a $0$-dimensional
component corresponding to a base point free $\mathfrak g^1_k$.
Moreover, if $\epsilon:\mm_{g-1,1}\rightarrow \mm_{g-1}$ is the
forgetful morphism, then using the identification
$\Delta_1=\mm_{g-1,1}\times \mm_{1,1}$, we have that
$\overline{Z}'\cap \Delta_1\supset \epsilon^*(Z)\times \mm_{1,1}$.
\end{theorem}
\begin{proof} We choose  a general
curve $[C]\in Z\subset GP_{g-1,k-1}^1$, a general point $p\in C$
and we set $X_0:=C\cup_p E$, where $E$ is an elliptic curve. By
assumption, there exists a base point free  $(A,V)\in
G^1_{k-1}(C)$ and $\rho\in \PP \mbox{Ker}(\mu_0(V))$ such that
$\mbox{dim}_{(A,V,\rho)} Q_{k-1}^1(C)=0$. In particular
$\mbox{Ker}(\mu_0(V))$ is $1$-dimensional and $h^0(A)=2$. Let
$\pi:X\rightarrow B$ be the versal deformation space of $X_0$,
$\Delta\subset B$ the boundary divisor corresponding to singular
curves, and we consider the scheme $\nu:\cQ_{k}^1\rightarrow B$
parametrizing $(gp)_{k}^1$-relations, which was constructed in
Theorem \ref{rel}. We construct a $(gp)_{k}^1$-relation $z=(l,m,
\rho_1, \rho_2)$ on $X_0$ as follows: the $C$-aspect of the limit
$\mathfrak g^1_k$ denoted by $l$ is obtained by adding $p$ as a
base point to $(A, V)$, while the $E$-aspect of $l$ is constructed
by adding $(k-2)p$ as a base locus to $|\OO_E(p+q)|$, where $q\in
E-\{p\}$ satisfies $2(p-q)\equiv 0$. Thus the vanishing sequences
$a^{l_C}(p)$ and $a^{l_E}(p)$  are $(1,2)$ and $(k-2,k-1)$
respectively. The $C$-aspect of the limit $\mathfrak
g_{2g-2-k}^{g-k}$ we denote by $m$, is the complete linear series
$|M_C|=|K_C(p)\otimes A^{-1}|$ which by Riemann-Roch has vanishing
sequence $(1,2, \ldots, g-k+1)$ at $p$. Finally the $E$-aspect of
$m$ is the subseries of $|\OO_E((2g-1-k)p-q)|$ with vanishing
$(g-3, g-2, \ldots, 2g-k-3)$ at $p$. From the base point free
pencil trick it follows that we can choose uniquely the relations
$\rho_C$ on $C$ and $\rho_E$ on $E$ such that
$\mbox{ord}_p(\rho_C)=3$ and $\mbox{ord}_p(\rho_E)=2g-5$ (we use
that $h^0(C, K_C\otimes A^{-2})=\mbox{dim}(\mbox{Ker}(\mu_0(V))=1$
by assumption, hence $\rho_C$ is essentially $\rho$ up to
subtracting the base locus).

 From Theorem \ref{rel}, every component of $\cQ_{k}^1$ passing
 through $z$ has dimension $\geq 3g-4$. On the other hand we claim that every
  component of $\nu^{-1}(\Delta)$ passing through $z$ has dimension
 $\leq 3g-5$ and that $z$ is an isolated point in $\nu^{-1}([X_0]
)$.
 Assuming this for a moment, we obtain that $z$ is a smoothable $(gp)_{k}^1$-relation in the sense
 that there is a component of $\cQ_{k}^1$ through $z$ which meets $\nu^{-1}(B-\Delta)$. From
 this it follows that $[X_0]\in \overline{GP}_{g,k}^1\cap
 \Delta_1$. Since by construction the curves $[X_0]$ fill up a divisor
 inside $\Delta_1$, we conclude that $GP_{g,k}^1$ has a
 divisorial component $Z'$ such that $\overline{Z}'\cap \Delta_1\supset \epsilon^*(Z)\times \mm_{1,1}$.

Furthermore, because the vanishing sequences of the $C$ and
$E$-aspects of $l$ add up precisely to $k$, every $\mathfrak
g^1_k$ on a  smooth curve \lq \lq near''$ X_0$  which specializes to $l$, is
base point free (cf. \cite{EH1}, Proposition 2.5). We obtain that a
point $z'\in \nu^{-1}(B-\Delta)$ near $z$ will satisfy
$\mbox{dim}_{z'}(\cQ_{k}^1)=3g-4$ and will correspond to a smooth
curve $[C']\in GP_{g,k}^1$, satisfying all the required
conditions.

We return now to the estimate for
$\mbox{dim}_z(\nu^{-1}(\Delta))$: we consider  a curve
$X_b=C_b\cup _q E_b$ with $b\in \Delta$, and let $(l,m,
\rho_{C_b}, \rho_{E_b})\in \nu^{-1}(b)$. Hence the underlying
limit linear series $l$ and $m$ have vanishing sequences
$a^{l_{C_b}}(q)=(1,2), a^{l_{E_b}}(q)=(k-2,k-1)$ and
$a^{m_{C_b}}(q)=(1,2,\ldots ,g-k+1),
a^{m_{E_b}}(q)=(g-3,g-2,\ldots,2g-3+k)$ respectively.
\newline
\indent Clearly $\mbox{ord}_q(\rho_{C_b})\geq 3(=1+2=2+1)$. We set
$\bigl(A_b=L_{C_b}(-q), V_{C_b}:=V_{C_b}(-q)\bigr)\in
G^1_{k-1}(C_b)$ and $\bigl(B_b=L_{E_b}(-(k-2)q),
V_{E_b}:=V_{E_b}(-(k-2)q\bigr)\in G^1_2(E_b)$. We claim that in
fact $\mbox{ord}_q(\rho_{C_b})=3$ and therefore
$\mbox{ord}_q(\rho_{E_b})=2g-5(=k-2+(2g-3-k)=k-1+(2g-4-k))$.
Indeed assuming that $\mbox{ord}_q(\rho_{C_b})\geq 4$, from the
base point free pencil trick we have that $h^0(C_b, K_{C_b}\otimes
A_b^{-2}(-q))\geq 1$. But $h^0(C,K_C\otimes A^{-2}(-p))=0$ (use
the assumption on $C$ and the fact that $p\in C$ is a general
point), which implies that we can assume that $h^0(C_b,
K_{C_b}\otimes A_b^{-2}(-q))=0$ for any point in a component of
$\nu^{-1}(\Delta)$ passing through $z$.

 After subtracting base points $\rho_{C_b}$ can be viewed as an element in the projectivization
of the kernel of the map $\mu_0(V_{C_b}):V_{C_b}\otimes
H^0(K_{C_b}\otimes A_b^{-1})\rightarrow H^0(K_{C_b})$, while
$\rho_{E_b}$ is in the projectivized kernel of the map
$\mu_0(V_{E_b}):H^0(E_b,B_b)\otimes H^0\bigl(E_b,
B_b^{-1}(4q)\bigr)\rightarrow H^0\bigl(E_b, \OO_{E_b}(4q)\bigr).$
In other words $[C_b]\in GP_{g-1,k-1}^1$ and from the base point
free pencil trick we get that $H^0(E_b, \OO_{E_b}(4q)\otimes
B_b^{-2})\neq 0$, which leaves only finitely many choices for
$B_b$ and $\rho_{E_b}$. It follows that $\mbox{dim}_z
\nu^{-1}(\Delta)\leq \mbox{dim}_{[C]}(GP_{g-1,k-1}^1)+1+1=3g-5$.
\end{proof}

\noindent \emph{Proof of Theorem \ref{petri}.} We apply Theorem
\ref{gp} starting with the base case $k\geq 3, g=2k-2$. In this
situation the locus $GP_{2k-2,k}^1$ is a divisor in $\cM_g$ which
can also be viewed as the branch locus of the map to $\cM_g$ from
the Hurwitz scheme of coverings $C\stackrel{k:1} \longrightarrow
\PP^1$ having a genus $g$ source curve (cf. \cite{EH3}, Section
5). The locus of $[C]\in \cM_g$ having infinitely many base point
free $\mathfrak g^1_k$'s is of codimension $\geq 2$, hence by
default the general point of $GP_{2k-2,k}^1$ corresponds to a
curve with finitely many $(A,V) \in G^1_k(C)$. The fact that for
each of these pencils, $\mbox{dim}\mbox{ Ker}(\mu_0(V))\leq 1$,
also follows from \cite{EH3}. Applying now Theorem \ref{gp}
repeatedly we construct divisorial components of
$GP_{2k-2+a,k+a}^1$ for all $k\geq 3$ and $a\geq 0$. It is easy to
check that in this way we fill all the cases claimed in the
statement. \hfill $\Box$

One could also define the loci $GP_{g,k}^1$ for $k\leq (g+1)/2$.
In this case $GP_{g,k}^1$ coincides with the locus of $k$-gonal
curves, which is irreducible of dimension $2g+2k-5$. When $g$ is
odd, $GP_{g, (g+1)/2}^1$ is the well-known Brill-Noether divisor
on $\cM_g$ introduced by Harris and Mumford (see \cite{EH3}). The
Gieseker-Petri divisors $GP_{g,k}^1$ with $k\geq (g+2)/2$ that we
introduced, share certain properties with the Brill-Noether
divisor. For instance the following holds (compare with
\cite{EH3}, Proposition 4.1):

\begin{proposition}
We denote by $j:\mm_{2,1}\rightarrow \mm_g$ the map obtained by
attaching a fixed general pointed curve $(C_0,p)$ of genus $g-2$.
Then for $(g+1)/2\leq k\leq g-1$ we have the relation
$j^*(\overline{GP}_{g,k}^1)=q \overline{\mathcal{W}}$, where
$q\geq 0$ and $\mathcal{W}$ is the divisor of Weierstrass points
on $\cM_{2,1}$.
\end{proposition}

\noindent \emph{Sketch of proof.}  We can degenerate $(C_0,p)$ to
a string of elliptic curves $(E_1\cup\ldots \cup E_{g-2}, p)$,
where $p$ lies on the last component $E_{g-2}$. We assume that for
all $2\leq i\leq g-2$, the points of attachment between $E_{i-1}$
and $E_i$ are general. Fix now $[B,p]\in \cM_{2,1}$ and assume
that $[X_0:=C_0\cup _p B] \in \overline{GP}_{g,k}^1$. We denote by
$(l_B, m_B, \rho_B)$ the $B$-aspect of a $(gp)_k^1$-relation on
$X_0$. Then using the setup described at the beginning of Section
4 we obtain that $\mbox{ord}_p(\rho_B)\geq 2g-4$. Since $l_B$ is a
$\mathfrak g^1_k$ and $m_B$ is a $\mathfrak g^{g-k}_{2g-2-k}$, the
only way this could happen is if $a^{l_B}(p)=(k-2,k)$ and
$a^{m_B}(p)=(\ldots,2g-4-k, 2g-2-k)$, which implies that
$h^0(\OO_B(2p))\geq 2$, that is, $[B,p]\in \mathcal{W}$. \hfill
$\Box$

\begin{remark}
Using  methods developed in  this section we can also prove
the following result  useful for the computation of the
class $[\overline{GP}_{g,k}^1]\in \mbox{Pic}(\overline{\cM}_g)$:
if $\epsilon:\mm_{g-1,1}\rightarrow \mm_{g-1}$ is the forgetful
morphism and $\phi:\mm_{g-1,1}\rightarrow \mm_{g}$ denotes the map
attaching an elliptic tail at the marked point, then
$\phi^*(\overline{GP}_{g,k}^1)$ is set-theoretically the union of
the two divisors: $\epsilon^*(\overline{GP}_{g-1,k-1}^1)$, and the
closure $D$ in $\mm_{g-1,1}$ of the locus of curves $[C,p]\in \cM_{g-1,1}$ for
which there exists a base point free $A\in W^1_k(C)$ such that
$h^0(C, A(-2p))\geq 1$ and the multiplication map $H^0(C,A)\otimes
H^0(K_C\otimes A^{\vee}(2p))\rightarrow H^0(K_C(2p))$ is not
injective. It is natural to view $D$ as a \lq \lq pointed"
Gieseker-Petri divisor on $\mm_{g,1}$.
\end{remark}

 We consider now the moduli space $\S_{g,n}$ of $n$-pointed
spin curves of genus $g$ and its subvariety $\S_{g,n}^r$
consisting of elements $(C, p_1,\ldots,p_n,L)$, where
$[C,p_1,\ldots, p_n]\in \cM_{g,n}$ and $L\in \mbox{Pic}^k(C)$ is a
line bundle such that $L^2\otimes \OO_C(p_1+\cdots +p_n)=K_C$ and
$h^0(L)\geq r+1$. Of course we assume that  $2k+n=2g-2$. The base
point free pencil trick relates these loci to the loci
$GP_{g,k}^i$ we introduced before. Precisely, if
$f:\S_{g,n}\rightarrow \cM_g$ is given by $[C, p_1,\ldots,
p_n,L]\mapsto [C]$, then $f(\S_{g,n}^1)= GP_{g,k}^1$.

We now look at the divisor $Z\subset GP_{g,k}^1$ constructed in
Theorem \ref{gp}. The condition that for a general $[C]\in Z$, the
scheme $Q_k^1(C)$ has a $0$-dimensional component with general
point corresponding to a base point free $\mathfrak g^1_k$, can be
translated into saying that $f^{-1}[C]$ has a zero-dimensional
component. We obtain in this way that there exists a component $Y$
of $\S_{g,n}^1$ of dimension $3g-4$ such that $f(Y)=Z$. This
proves Corollary \ref{pol}.

\section{Injectivity of Gaussian maps}

We are going to prove Theorem \ref{th3} by degeneration. Our proof
is inspired by the work of Eisenbud and Harris on the
Gieseker-Petri Theorem (cf. \cite{EH2}). Suppose we have a family
of genus $g$ curves $\pi:X\rightarrow B$ over a base
$B=\mbox{Spec}(R)$ with $R$ being a complete DVR with local
parameter $t$ and let $0$ and $\eta$ respectively, denote the
special and the generic point of $R$. Assume furthermore that
$X_{\eta}$ is smooth and that $X_0$ is a curve of compact type
consisting of a string of components of which $g$ of them $E_1,
\ldots E_g$, are elliptic curves, while the rest are rational
curves, glued in such a way that the stable model of $X_0$ is the
curve $E_1\cup _{p_1} E_2\cup_{p_2} E_3\cup \ldots \cup
E_{g-1}\cup_{p_{g-1}} E_g.$ Slightly abusing the notation, for
$2\leq i\leq g-1$ we will consider  $p_{i-1}$ and $p_i\in E_i$ to
be the points of attachment of $E_i$ to $\overline{X_0-E_i}$ and
we will choose $X_0$ in such a way that $p_{i}-p_{i-1}$ is not a
torsion class in $\mbox{Pic}^0(E_i)$.

We proceed by contradiction and assume that there exists a line
bundle $L_{\eta}$ on $X_{\eta}$ of degree $d$, together with a
non-zero element $$\rho_{\eta} \in \mbox{Ker}\{\psi_{L\eta}:\wedge
^2H^0(X_\eta, L_\eta)\rightarrow H^0(X_\eta,
\Omega^1_{X_{\eta}}\otimes L_{\eta}^2)\}.$$ (Note that we because
the shape of $X_0$ does not change if we blow-up the surface $X$,
we can assume that we have a bundle $L_{\eta}$ on $X_{\eta}$
rather than on the geometric generic fibre $X_{\overline{\eta}}$.)
As in Section 4, for each component $Y$ of $X_0$ we have the line
bundle $L^Y$ on $X$ extending $L_{\eta}$ and having degree $0$
restriction to all components $Z\neq Y$ of $X_0$ and we
 set $L_Y:=L^Y_{|Y}$.
Starting with $\rho_{\eta} \in \wedge^2 \pi_*(L_{\eta})$  we
obtain elements $\rho^Y =t^{\alpha}\rho_{\eta} \in \wedge ^2
\pi_*(L^Y)-t\wedge ^2\pi_*(L^Y)$ for uniquely determined integers $\alpha$, and we define
$\rho_Y:=\rho^Y_{|Y} \in \wedge ^2 H^0(Y, L_Y)$.
\begin{lemma}\label{y}
For each component $Y$ of $X_0$ we have that \begin{center}
$\rho_Y \in
\rm{Ker}$$\{\psi_{L_Y}:\wedge^2 H^0(Y,L_Y)\rightarrow H^0(Y,
\Omega_Y^1 \otimes L_Y^2)\}.$
\end{center}
\end{lemma}
\begin{proof} We use the commutative diagram
$$\begin{array}{ccc} \wedge^2 H^0(X_0, L^Y_{|X_0}) &
\stackrel{res}\longrightarrow & \wedge^2 H^0(Y,
L_Y) \\
\rmapdown{\psi_{L^Y_{|X_0}}} & \; & \rmapdown{\psi_{L_Y}}\\
H^0(X_0,\Omega_{X_0}^1\otimes L_{|X_0}^{Y^{\otimes ^2}}) &
\stackrel{res}\longrightarrow &
H^0(Y,\Omega_Y^1\otimes L_Y^2) \\
\end{array}$$
and keep in mind that the upper restriction map is injective.
\end{proof}

We will use the following observation (similar to the one for
ordinary multiplication maps): let $C$ be a smooth curve, $p\in C$
and $M$ a line bundle on $C$. If $\rho \in \mbox{Ker}(\psi_M)$ and
$\{\sigma_i\}$ is a basis of $H^0(M)$ such that
$\mbox{ord}_p(\sigma_i)=a_i^M(p)=a_i$, then there are distinct
pairs of integers $(i_1,j_1)\neq (i_2,j_2)$ with $i_1\neq j_1$ and
$i_2\neq j_2$, such that
$\mbox{ord}_p(\rho)=\mbox{ord}_p(\sigma_{i_1})+\mbox{ord}_p(\sigma_{j_1})=\mbox{ord}_p(\sigma
_{i_2})+\mbox{ord}_p(\sigma_{j_2})$. This follows from a local
calculation: if $t$ is a local parameter for $C$ at $p$, then
$$\psi_M(\sigma_i \wedge \sigma_j)=\bigl((a_i-a_j)
t^{a_i+a_j-1}+\mbox{ h.o.t.}\bigr)dt, $$
and since $\psi_M(\rho)=0$, the number $\mbox{ord}_p(\rho)$ must be attained
for at least two pairs $(i,j)$.

\begin{proposition}\label{ind}
Suppose $Y$ and $Z$ are two components of $X_0$ meeting at a point
$q$ and let $p$ be a general point on $Y$. We have the following
inequalities:
\begin{enumerate}
\item $\rm{ord}$$_q(\rho_Z)\geq \rm{ord}$$_p(\rho_Y).$ \item If
$Y$ is one of the elliptic components of $X_0$, then
$\rm{ord}$$_q(\rho_Z)\geq \rm{ord}$$_p(\rho_Y)+2.$
\end{enumerate}
\end{proposition}
\begin{proof}
Although (1) is essentially Proposition 3.1 from \cite{EH2} we
will briefly go through the proof and in doing so we will also prove (2). We
pick a basis $\{\sigma_i=\sigma_i^Y\}$ of $\pi_*(L^Y)$ such that
$\mbox{ord}_p(\sigma_{i |Y})=a_i^L(p) $ and for which there are
integers $\alpha_i$ with the property that
$\{\sigma_i^Z=t^{\alpha_i}\sigma_i\}$ form a basis for $H^0(X,
L^Z)$ (see \cite{EH1}, Lemma 2.3 for the fact that such a basis
can be chosen). We then write $\rho^Y=\sum_{i\neq j} f_{ij}
\sigma_i\wedge \sigma_j$, with $f_{ij}\in R$, and we can express
$\rho^Z=t^{\gamma}\rho^Y$, where
$\gamma=\mbox{max}_{i\neq j}\{\alpha_i+\alpha_j-\nu(f_{ij})\} $. Here
$\nu$ denotes the valuation on the ring $R$. From the definition
of $\gamma$ it follows that there exists a pair $(i,j), i\neq j$, with
$\gamma=\alpha_i+\alpha_j-\nu(f_{ij})$, such that we have a string
of inequalities
\begin{equation}\label{e2}
\mbox{ord}_q(\rho_Z)=\mbox{ord}_q(\sigma_{i |Z}^Z)+\mbox{ord}_q(\sigma_{j |Z}^Z)\geq
\alpha_i+\alpha_j\geq \gamma, \end{equation}
(see also Section 4). On the other hand
there exists a  pair $(i',j'), i'\neq j'$ such that $\nu(f_{i'j'})=0$, for
which we can write the inequalities
\begin{equation}\label{e3}
 \mbox{ord}_p(\rho_Y)=
 \mbox{ord}_p(\sigma_{i'| Y})+\mbox{ord}_p(\sigma_{j'| Y})\leq
(d-\mbox{ord}_q(\sigma_{i' | Y}))+(d-\mbox{ord}_q(\sigma_{j' | Y}))\leq
\alpha_{i'}+\alpha_{j'}\leq \gamma. \end{equation} Combining
(\ref{e2}) and (\ref{e3}) we get the first part of the
Proposition. When moreover the curve $Y$ is elliptic, since
$\psi_{L_Y}(\rho_Y)=0$, there must at least two pairs $(i_1,j_1)$
and $(i_2,j_2)$ for which (\ref{e3}) holds. On the other hand
 $p-q \in \mbox{Pic}^0(Y)$ can be assumed not to be a torsion
class,  and we obtain that $\mbox{ord}_p(\sigma_{i |
Y})+\mbox{ord}_q(\sigma_{i | Y})\leq d-1$ for all indices $i$
except at most one. This and the fact that the vanishing orders
$\mbox{ord}_p(\sigma_{i | Y})$ are all distinct, quickly lead to
the inequality $\mbox{ord}_q(\rho_Z)\geq \gamma \geq
\mbox{ord}_p(\rho_Y)+2.$
\end{proof}
A repeated application of Proposition \ref{ind} gives the
following result:
\begin{proposition}\label{central}
Let $X_0$ be the curve described in the degeneration above and
which has the stable model $\cup_{i=1}^g E_i$, where $E_i$ are
elliptic curves. We denote by $p_{i-1}$ and $p_i$ the points of
attachment of $E_i$ to the rest of $X_0$.  If
$\psi_{L_{\eta}}(\rho_{\eta})=0$, then
$\rm{ord}$$_{p_{g-1}}(\rho_{E_g})\geq
\rm{ord}$$_{p_1}(\rho_{E_2})+2g-4.$
\end{proposition}
We are now in a position to prove Theorem \ref{th3}. In fact we have a more
general result:
\begin{theorem}\label{th4}
For a general genus $g$ curve $C$ and for any line bundle $L$ on $C$ of
degree $d\leq a+g+2$, where $a\geq 0$, we have
that $\rm{dim }$ $\rm{Ker}$$(\psi_L)\leq a(a+1)$. In particular,
if $d\leq g+2$ then $\psi_L$ is injective.
\end{theorem}
\begin{proof} We apply Proposition \ref{central} and degenerate $C$ to $X_0=E_1\cup \ldots \cup E_g$. We assume that
$\mbox{Ker}(\psi_{C,L})$ is at least $1+a(a+1)$-dimensional. Then
  $\mbox{dim } \mbox{Ker}(\psi_{X_0,L^{E_2}_{\ |X_0}})\geq
1+a(a+1)$  and since the restriction map $\wedge^2
H^0(X_0,L^{E_2}_{\ |X_0})\rightarrow \wedge ^2 H^0(E_2,L_{E_2})$
is injective we  obtain that $\mbox{Ker}(\psi_{E_2,L_{E_2}})$ is
at least $1+a(a+1)$-dimensional as well. For simplicity let us
denote $E_2=E, L_{E_2}=L$ and $p_1=p\in E_2$ (recall that $p_1\in
E_2\cap E_1).$

If we choose a basis $\{\sigma_i\}$ of $H^0(L)$ adapted to the
point $p$, then as we noticed before for each $\rho \in
\mbox{Ker}(\psi_L)$ there will be at least two distinct pairs of
integers $(i_1,j_1)\neq (i_2,j_2)$ where $i_1\neq j_1, i_2\neq
j_2$ such that
$$\mbox{ord}_p(\rho)=\mbox{ord}_p(\sigma_{i_1})+\mbox{ord}_p(\sigma_{j_1})
=\mbox{ord}_p(\sigma_{i_2})+\mbox{ord}_p(\sigma_{j_2}).$$
The vanishing sequence  $a^{L_{E_1}}(p)$ is $\leq
(\ldots, d-3,d-2,d)$, hence  the vanishing sequence
of $L=L_{E_2}$ at $p$ is $\geq (0,2,3,4,5,\ldots)$, which yields
that $\mbox{ord}_p(\rho)\geq 5(=0+5=2+3)$ for every $\rho \in
\mbox{Ker}(\psi_L)$. Since $\mbox{dim }\mbox{Ker}(\psi_L)\geq 1+a(a+1)$,  there is a subspace $W_1\subset
\mbox{Ker}(\psi_L)$ of dimension $\geq a(a+1)$ such that
$\mbox{ord}_p(\rho)\geq 6(=0+6=2+4)$ for each $\rho \in W_1$.
\newline
\indent Repeating this reasoning for $W_1$ instead of
$\mbox{Ker}(\psi_L)$ we obtain a subspace $W_2\subset W_1$ with
$\mbox{dim}(W_2)\geq \mbox{dim}(W_1)-1$ such that
$\mbox{ord}_p(\rho)\geq 7(=0+7=2+5=3+4)$ for every $\rho \in W_2$,
and then a subspace $W_3\subset W_2$ with $\mbox{dim}(W_3)\geq
\mbox{dim}(W_2)-2$ with the property that $\mbox{ord}_p(\rho)\geq
8(=0+8=2+6=3+5)$ for all $\rho \in W_3$. At the end of this argument we find at
least one element $\rho=\rho_{E_2}\in \mbox{Ker}(\psi_L)$ such
that $\mbox{ord}_p(\rho)\geq 2a+5$. Since this reasoning works if
we replace $\mbox{Ker}(\psi_L)$ with any of its subspaces having
dimension $\geq 1+a(a+1)$, we can assume that $\rho_{E_2}$ is the
restriction to $E_2$ of an element $\rho_{\eta}$ in the kernel of the
corresponding Gaussian map on the general curve $X_{\eta}$, which
according to the procedure described before Lemma \ref{y} will
produce elements $\rho_{E_i}\in \mbox{Ker}(\psi_{L_{E_i}})$ for
$1\leq i\leq g$. Applying Proposition \ref{central} we have that
$\mbox{ord}_{p_{g-1}}(\rho_{E_g})\geq
\mbox{ord}_{p}(\rho)+2g-4=2(a+g)+1$. The vanishing sequence of
$L_{E_g}$ at $p_g$ is $\leq (\ldots,d-3,d-2,d)$ from which we
obtain that on the other hand
$$\mbox{ord}_{p_{g-1}}(\rho_{E_g})\leq 2d-5(=d+(d-5)=(d-2)+(d-3)),$$
which combined with the previous inequality yields $d\geq a+g+3$
which is a contradiction.
\end{proof}

Note that Theorem \ref{th4} is valid for an \emph{arbitrary} line
bundle on a general genus $g$ curve. It is clear that Proposition
\ref{central} would give better sufficient conditions for the
injectivity of $\psi_L$ if we restricted ourselves to line bundles
on $C$ having a prescribed ramification sequence at a given point
$p\in C$. In this case we degenerate $(C,p)$ to $(X_0=E_1\cup
\ldots \cup E_g,p)$, where $X_0$ is as in Theorem \ref{th4} and
$p\in E_1$ is such that $p-p_1\in \mbox{Pic}^0(E_1)$ is not a
torsion class. We leave it to the interested reader to work out
the numerical details. We can also improve on Theorem \ref{th4} if
we look only at a suitably general line bundle $L$ on $C$:

\begin{proposition}
Fix integers $g,d$ and $r\geq 2$ such that $d\leq g+r$, $\rho=g-(r+1)(g-d+r)\geq 0$ and
moreover $d<g+3+\frac{\rho}{2(r-1)}.$ Then if $C$ is a general
curve of genus $g$ and $L\in W_d^r(C)$ is general,  the Gaussian map
$\psi_L$ is injective.
\end{proposition}
\begin{proof} We degenerate $C$ to $X_0$, fix a general point $p\in
E_1$ and set $a:=[\rho/(r-1)]+2$. Our numerical assumptions imply
that $\rho- (a-2)(r-1)\geq 0$. From the general theory of limit
linear series in \cite{EH1} reducing the Brill-Noether theory of
$X_0$ to Schubert calculus, we know that there exists a smoothable limit
linear series of type $\mathfrak g^r_d$ on $X_0$, say
$l=\{L_{E_i}\in W^r_d(E_i)\}_{i=0,\ldots,g}$ having vanishing
sequence  $\geq (0,1, a, a+1, a+2, \ldots,a+r-2)$ at the point $p$.

Assume by contradiction that there are elements $\rho_{E_i}\in
\mbox{Ker}(\psi_{L_{E_i}})$ coming from an element  $\rho \neq 0$
in the kernel of the corresponding Gaussian on the general curve.
Then $\mbox{ord}_p(\rho_{E_1})\geq a+1(=1+a=0+(a+1))$ and from
Proposition \ref{ind} we get that
$\mbox{ord}_{p_{g-1}}(\rho_{E_g})\geq
\mbox{ord}_p(\rho_{E_1})+2g-2=2g+a-1$. On the other hand, as we noticed
before $\mbox{ord}_{p_{g-1}}(\rho_{E_g})\leq 2d-5$ which gives a
contradiction.
\end{proof}
\begin{remark} The techniques from this section also allow us to study the
kernel $S_2(L)$ of the multiplication map
$\mu_L:\mbox{Sym}^2H^0(L)\rightarrow H^0(L^2)$. In a way similar
to the proof of Theorem \ref{th4} we can show that if $L$ is an
arbitrary line bundle of degree $d\leq g+a+1$ on a general curve
$C$ of genus $g$ then $\mbox{dim }S_2(L)\leq a(a+1)$. The $a=0$
case of this result has been established by Teixidor (cf.
\cite{T2}). We also note that this result as well as Theorem
\ref{th4}, are meaningful when the bundle $L$ is special. On the
other hand the case when $L$ is nonspecial (when, under suitable
assumptions, we expect surjectivity for both $\psi_L$ and
$\mu_L$), has been extensively covered in the literature (see e.g.
\cite{Pa}).
\end{remark}
 Theorem \ref{th3} answers Question 5.8.1 from Wahl's survey
\cite{W1}, where the problem is raised in terms of
self-correspondences on a curve. Suppose that $C$ is a smooth
curve and we consider the diagonal $\Delta \subset C\times C $ and
the projections $p_i:C\times C\rightarrow C$ for $i=1,2$. For a
line bundle $L$ on $C$ we denote $L_i:=p_i^*L$ for $i=1,2$. We can
rephrase Theorem \ref{th3} as follows:
\begin{proposition}\label{corr}
If $L$ is a line bundle of degree $d\leq g+1$ on a general curve
$C$ of genus $g$, then $H^0(C\times C, L_1+L_2-2\Delta)=0$.
\end{proposition}
\begin{proof} We use that $H^0(C\times C,
L_1+L_2-2\Delta)=\mbox{Ker}(\Phi_L)=S_2(L)\oplus
\mbox{Ker}(\psi_L)$. We have proved that $\mbox{Ker}(\psi_L)=0$
while $S_2(L)=0$ follows from \cite {T2}.
\end{proof}

\end{document}